\documentclass[a4, 11pt]{article}
\usepackage{amssymb}
\usepackage{amsmath}
\usepackage{bbm}
\usepackage{empheq}
\usepackage[framemethod=tikz]{mdframed}
\usepackage{lipsum}
\usepackage{tcolorbox}
\usepackage{fourier-orns}
\usepackage{enumitem}
\RequirePackage[amsmath,thmmarks,hyperref,framed]{ntheorem}
\usepackage[T1]{fontenc}
\usepackage{lmodern}

\usepackage[applemac]{inputenc}

\newtheorem{theorem}{Theorem}
\newtheorem{prop}[theorem]{Proposition}
\newtheorem{cor}[theorem]{Corollary}
\newtheorem{defi}[theorem]{Definition}
\newtheorem{lemma}[theorem]{Lemma}

\numberwithin{equation}{section}

\usepackage[english]{babel}
\newmdenv[innerlinewidth=0.5pt, roundcorner=4pt,innerleftmargin=6pt,
innerrightmargin=6pt,innertopmargin=10pt,innerbottommargin=10pt,backgroundcolor=gray!21]{mybox}

\newcommand{\abs}{|\alpha|}

\newcommand{\T}{T_{\mathrm{sep}}}
\newcommand{\alphaF}{\emph{\texttt{a}}}

\theorembodyfont{\upshape}

\date{}
\title{Tail asymptotics for extinction times of self-similar fragmentations}

\author{B\'en\'edicte Haas\thanks{Universit\'e Sorbonne Paris Nord, LAGA, CNRS (UMR  7539) 93430 Villetaneuse, France \newline \hspace*{0.5cm} E-mail: haas@math.univ-paris13.fr}}

\begin{document}

\maketitle

\begin{abstract}
We provide the exact large-time behavior of the tail distribution of the extinction time of a self-similar fragmentation process with a negative index of self-similarity,  improving thus a previous result on the logarithmic asymptotic behavior of this tail. Two factors influence this behavior: the distribution of the largest fragment at the time of a dislocation and the index of self-similarity.  As an application we obtain the asymptotic behavior of all moments of the largest fragment and compare it to the behavior of the moments of a tagged fragment, whose decrease is in general significantly slower. We illustrate our results on several examples, including fragmentations related to random real trees -- for which we thus obtain the large-time behavior of the tail distribution of the height -- such as the stable L\'evy trees of Duquesne, Le Gall and Le Jan (including the Brownian tree of Aldous), the alpha-model of Ford and the beta-splitting model of Aldous.
\end{abstract}

%%%%%%%%%%%%%%%%%%%%%%
\section{Introduction}
%%%%%%%%%%%%%%%%%%%%%%

The aim of this article is to give exact bounds for the (upper) tail distribution of the extinction time of a self-similar fragmentation with a negative index of self-similarity.  In Probability theory, self-similar fragmentation processes, which are meant to model the evolution of a system of particles that fragment repeatedly, were developed in a general framework including infinite splitting rates by Bertoin at the beginning of the 2000's \cite{BertoinHomogeneous,BertoinSSF,BertoinAB, BertoinBook}, inspired by previous works on fragmentation and coalescence processes (to name a few: \cite{Kolmo41,Filippov, BD86,BD87,King82,PitCoag99}). In these models of fragmentations with a \emph{negative} index of self-similarity, particles are subject to successive fragmentations which intensify as the masses of particles decrease. This leads to a shattering phenomenon \cite{McGZ87, Filippov, BertoinAB}, where a particle is entirely reduced to dust --  a set of zero-mass particles -- in finite time which is what we call the extinction time of the particle. Bounds for the \emph{logarithm} of the tail of the distribution of this extinction time where established in \cite{HaasLossMass}, under a regular variation assumption on the parameters of the model. Our objective here is to refine considerably this result by setting up precise bounds for the \emph{non-logarithmic} tail of the extinction time, within a broad framework of parameters, relaxing the regular variation assumption. As a byproduct we will obtain the large time behavior of all moments of the largest fragment. Interestingly their decrease is significantly slower (this is explicitly quantified) than that of moments of a typical fragment, which was already known. These results also allow us to conjecture the existence of a Yaglom limit for a rescaled version of the fragmentation process conditioned on non extinction.

\bigskip

In practice, a difficulty is to concretize on specific examples the overall rather abstract results that we obtain. We will illustrate a fairly general approach on a series of `natural' examples, for most of which the results we get are new.  This unified framework also allows us to retrieve estimates already known (\cite{DW17}) for the tail distributions of the heights of stable L\'evy trees.

\bigskip

\textbf{Organization of the paper.} In Section \ref{sec:intro} we gather background on self-similar fragmentations and set the notation that we will use throughout the paper. We then state in Section \ref{sec:mainresults} our main results on the asymptotic behavior of the tail of the extinction time of a fragmentation with a negative index of self-similarity, and the moments of the largest fragment. The proofs of these results are undertaken in Section \ref{sec:ProofTh1}, Section \ref{sec:coro} and Section \ref{sec:proofdecrease}. Section \ref{sec:examples} is then devoted to a series of illustrative examples. 

\bigskip

\textbf{Notation.} For  $f,g:\mathbb R_+ \rightarrow \mathbb R_+$ two non negative functions, we write
$$
f  ~ \asymp  ~ g \qquad (\text{ or } f (t)  ~ \asymp  ~ g (t))
$$
if there exist $c_1,c_2>0$ such that $c_1 g(t)\leq f(t) \leq c_2 g(t)$ for all $t$ large enough. If only one of the inequality holds, we write
$$
f  ~\lesssim  ~g  \quad \text{or}  \quad f  ~ \gtrsim  ~ g
$$
according to whether $f(t) \leq c g(t)$ for some $c>0$ and all $t$ large enough, or $f(t) \geq c g(t)$ for some $c>0$ and all $t$ large enough. Last if the functions are asymptotically proportional, i.e. if the ratio $f(t)/g(t) \rightarrow c$ for some $c>0$ as $t \rightarrow \infty$, we write
$$
f ~ \propto ~ g.
$$

%%%%%%%%%%%%%%%%%%%%%%%%%%%%%%%%%%%%%%%
\section{Background on self-similar fragmentations}
\label{sec:intro}
%%%%%%%%%%%%%%%%%%%%%%%%%%%%%%%%%%%%%%%

\textbf{($\mathbf{\alpha,\nu}$)-fragmentations.} Different state spaces have been considered to model self-similar fragmentations, depending on whether one wishes to keep track of the genealogy of the fragmentation or is interested only in the evolution of the masses of the particles involved. We refer to the book of Bertoin \cite{BertoinBook} for a detailed introduction to these processes.  Here we will mainly be interested in the evolution of masses, and so we will work on the set of non-increasing non-negative sequences with sum less or equal to 1:
$$
\mathcal S^{\downarrow}:=\left\{\mathbf s=(s_i)_{i \geq 1}\in [0,1]^{\infty} : \sum_{i=1}^{\infty} s_i \leq 1 \right\}.
$$
The formal definition of a self-similar fragmentation with index $\alpha \in \mathbb R$ then holds as follows:
 
\begin{defi} A \emph{self-similar fragmentation process $F$ with index $\alpha \in \mathbb R$} is a $\mathcal S^{\downarrow}$-valued Markov process continuous in probability such that if $\mathbb P_{\mathbf s}$ denotes the distribution of $F$ starting from $\mathbf s =(s_i)_{i \geq 1 }\in \mathcal S^{\downarrow}$, then
$$
F \text{ under } \mathbb P_{\mathbf s} ~  \overset{(\mathrm{d})}  = ~ \left(\left \{ s_iF^{(i)}(s_i^{\alpha}t)\right\}^{\downarrow}, t \geq 0\right)
$$
where the $F^{(i)}, i\geq 1$ are i.i.d. with common distribution $\mathbb P_{(1,0,\dots)}$ and the descending arrow means that the elements involved in the brackets are ranked in decreasing order. 
\end{defi}

The dynamics of such a process is therefore characterized by a branching property (given their masses, different particles evolve independently) and a scaling property (a particle fragments at a rate proportional to its mass to the power $\alpha$). In the following we will group these two properties under the name of \emph{fragmentation property}.

\medskip

Throughout we will focus on self-similar fragmentation processes that are \emph{pure-jump} (i.e. with no erosion, with Bertoin's vocabulary) and with a \emph{negative index of self-similarity} $\alpha$. Also, we will always work under $\mathbb P := \mathbb P_{(1,0,\ldots)}$. The distribution of such a process is then characterized by two parameters \cite{BertoinSSF, Berest02}: the index $\alpha<0$ and a $\sigma$-finite measure $\nu$ on  $\mathcal S^{\downarrow}$ such that $\int_{\mathcal S^{\downarrow}}(1-s_1)\nu(\mathrm d \mathbf s)<\infty$, generally called the \emph{dislocation measure}. 
We can summarize roughly the roles of $\alpha$ and $\nu$ as follows: 
\begin{center}
a particle of mass $m$ fragments at rate $m^{\alpha}\nu(\mathrm d \mathbf s)$.
\end{center}
When the measure $\nu$ is finite this means that a particle of mass $m$ waits an exponential time with parameter $m^{\alpha}$ before splitting into sub-particles with masses $m\mathbf s$, with $\mathbf s$ distributed according to $\nu(\cdot )/\nu(\mathcal S^{\downarrow})$. For more general measures $\nu$, such a process can be constructed via a Poisson point process and appropriate time changes depending on the ancestral lines of the fragments, see the above mentioned references. 

\medskip

\emph{From now on $F$ denotes such an $(\alpha,\nu)$-fragmentation process, with $\alpha<0$ and $F(0)=(1,0,\ldots)$}. We also assume that $\nu(\mathcal S^{\downarrow})>0$ so that the particles indeed fragment. All these assumptions are implicit is the rest of the paper.

\bigskip

\textbf{The extinction time $\zeta$.} 
As mentioned in the Introduction, since $\alpha<0$, the fragmentation of small particles intensifies as time goes on and reduces to dust the whole initial mass after a finite time, see in particular \cite[Proposition 2 (i)]{BertoinAB} for a proof in the general framework of self-similar fragmentations considered here, and \cite{McGZ87, Filippov} for previous references. We denote by 
\begin{equation*}
\label{def_zeta}
\zeta:=\left \{t \geq 0 :F(t)=0 \right\}
\end{equation*}
the first time at which the process is entirely reduced to dust. Our goal is to obtain a precise estimate of ~$\mathbb P(\zeta > t)$ as $t \rightarrow \infty$ in a fairly general framework.

\bigskip

\textbf{The functions $\phi$ and $\psi$, and connections with a tagged fragment.} The asymptotic behavior of the tail $\mathbb P(\zeta>t)$ will depend on a function $\psi$ which is built from the measure $\nu$ and that we now introduce. First we introduce the function $\phi$
\begin{equation}
\label{def:phi}
\phi(x):=\int_{\mathcal S^{\downarrow}} \bigg(1-\sum_{i \geq 1} s_i^{x+1} \bigg) \nu (\mathrm d \mathbf s), \quad x \geq 0
\end{equation}
which is the Laplace exponent of a subordinator involved in the distribution of the evolution of a typical fragment of the process. To settle this formally, we recall that it is possible to associate to the process $F$ an interval version of this fragmentation, i.e. a process $(O_F(t), t \geq 0)$
of nested open subsets of $]0,1[$ ($O_F(t_2) \subset O_F(t_1)$ for all times $t_2 \geq t_1 \geq 0$) such that $F(t)$ is the decreasing rearrangement of the lengths of the interval components of $O_F(t)$, for all $t\geq 0$, see the papers \cite{BertoinSSF,Berest02}. In such an interval version, consider a tagged point $U$ uniformly distributed on $(0,1)$ and independent of $O_F$ and let $F_{\mathrm{tag}}(t)$ denote the size of the fragment  (i.e. the interval component) containing $U$ at time $t$. Then the process $F_{\mathrm{tag}}$ writes
\begin{equation}
\label{lambda_sub}
F_{\mathrm{tag}}(t)=\exp(-\xi_{\rho(t)}), \quad t \geq 0
\end{equation}
where $\xi$ is a subordinator with Laplace transform $\phi$ and $\rho$ the Lamperti time-change defined for all $t\geq 0$ by $\rho(t)=\inf{\{u \geq 0:\int_0^u \exp(-\abs \xi_r) \mathrm dr >t\}}$, see  \cite[Chapter 3.2]{BertoinBook} or \cite[Section 4]{BertoinSSF}. Note that the process $F_{\mathrm{tag}}$ reaches 0 at time 
\begin{equation}
\label{def_I}
\zeta_{\mathrm{tag}}:=\int_0^{\infty} \exp(-\abs \xi_r) \mathrm dr,
\end{equation}
which belongs to the well-studied family of \emph{exponential functionals of subordinators}. This random variable is called the extinction time of the tagged fragment. Note also, for further calculations, that for each fixed $t$, conditional on $F(t)$, $F_{\mathrm{tag}}(t)=F_i(t)$ with probability $F_i(t)$ for all $i \geq 1$, and $F_{\mathrm{tag}}(t)=0$ with probability $1-\sum_{i\geq 1} F_i(t)$.

\bigskip

The function $\psi$ we are interested in is now defined as the inverse of the increasing function $x \mapsto x/\phi(x)$, which is a bijection from $(0,\infty)$ to $(x_{\psi}, \infty)$, where $x_{\psi}:=\big(\int_{\mathcal S^{\downarrow}}\sum_{i\geq 1}s_i\ln(s_i) \nu(\mathrm d \mathbf s))\big)^{-1}$ is the inverse of the right derivative at 0 of $\phi$. So for all $x>x_{\psi}$
\begin{equation}
\label{def:psi}
\frac{\psi(x)}{\phi(\psi(x))}=x.
\end{equation}
Note that the functions $\psi$ and $x \mapsto \psi(x)/x$ are both increasing since, as a Laplace exponent of a subordinator, $\phi$ is increasing and $x \mapsto \phi(x)/x$ is decreasing. 

\bigskip

A point that will be crucial in our study is the asymptotic behavior of the tail $\mathbb P(\zeta_{\mathrm{tag}}>t)$. As said above, exponential integrals of subordinators have been well-studied, see  \cite{BYsurvey} for a survey and \cite{MZ06,R12, HR12,MS21, H20a} for results on the large-time behavior of the tail or density of such random variables. In particular, the recent papers  \cite{MS21, H20a} give a precise estimate, assuming the forthcoming hypothesis (\ref{hyp:main}) (see next section) on the Laplace exponent of the subordinator. With our notation and focusing on the random variable $\zeta_{\mathrm{tag}}$, this estimates reads
\begin{equation}
\label{eq:tailI}
\mathbb P(\zeta_{\mathrm{tag}}>t) \  \propto \  \frac{t (\psi'(\abs t))^{1/2}}{\psi(\abs t)} \exp\left(- \int_{x_{\psi}+1}^{t} \frac{\psi(\abs r)}{\abs r} \mathrm dr\right)
\end{equation}
and we have that $\zeta_{\mathrm{tag}}$ is in the domain of attraction of an exponential distribution:
\begin{equation}
\label{extreme}
\frac{\mathbb P\left(\zeta_{\mathrm{tag}} >t+x\frac{\abs t}{\psi(\abs t)}\right)}{\mathbb P(\zeta_{\mathrm{tag}}>t)} ~ \underset{t \rightarrow \infty}\longrightarrow ~ \exp(-x), \quad \forall x \geq 0.
\end{equation}
This last point is certainly a consequence of  (\ref{eq:tailI}), but it was actually known previously and is closely related to the large-time behavior of $F_{\mathrm{tag}}$ conditioned on being positve. See \cite[Proposition 2.7 and Proposition 3.1]{HR12} for details.

\bigskip

\textit{Remark.} Note that if we set $\phi_1(x):=\int_{\mathcal S^{\downarrow}}(1-s_1^{x+1})\nu(\mathrm d \mathbf s)$, then as $x \rightarrow \infty$
$$
\phi(x)=\phi_1(x)+O\left(\left(\frac{1}{2}\right)^{x}\right) \quad \text{ and } \quad \phi^{(k)}(x)= \phi^{(k)}_1(x)+ O\bigg(\left(\frac{1}{2}\right)^{x-\varepsilon}\bigg)
$$
for the derivatives at any order $k \geq 1$ and for all $\varepsilon>0$. Consequently, defining $\psi_1$ as the inverse of the function $x \mapsto x/\phi_1(x)$ (which is well-defined in a neighborhood of $+\infty$) we can replace the function $\psi$ by the function $\psi_1$ in the above statement (\ref{eq:tailI}) and in the forthcoming Theorem \ref{thm:main}, Proposition \ref{prop-eps-2} and all their corollaries (including Proposition \ref{prop:decreaseF1}).

\bigskip

\textbf{Bounds for the logarithm of $\mathbb P(\zeta>t)$.} With the notation introduced above and under the assumption that $\phi$ is nearly regularly varying at $\infty$ (in the sense that $\phi \asymp f$, with $f$ regularly varying), it was proved in \cite[Prop 14]{HaasLossMass} that
$$
\ln \mathbb P(\zeta>t) ~\asymp ~-\psi(t).
$$
In particular the tail of $\zeta$ is exponential or even lighter (in fact this fast decrease holds without any hypothesis on $\phi$). 

\medskip

For comparison, we mention that the tail of the extinction time of a non-monotonic \emph{self-similar growth-fragmentation process} -- which is meant to describe the evolution of particles whose masses may vary `between' fragmentation times -- typically decreases as a power function (see \cite[Corollary 4.5]{BBCK18} and the references therein for background on growth-fragmentation processes). We note that the study in this case is quite different from that in the monotonic case and that the authors obtain quite easily a non-logarithmic estimate of the tail.  

\bigskip

\textbf{First examples of non-logarithmic estimates.} For the self-similar fragmentations without growth that we study here, there are only a few special cases in the litterature where a non-logarithmic asymptotic expression of $\mathbb P(\zeta>t)$ in terms of simple functions is available. In particular, there is the important \emph{Brownian fragmentation} (introduced in \cite[Section 4]{BertoinSSF}): in this case the extinction time $\zeta_{\mathrm{Br}}$ is distributed as the maximum of a Brownian excursion of length 1 and it therefore holds (see \cite{Kennedy}) that
$$
\mathbb P(\zeta_{\mathrm{Br}}> t) \underset{t \rightarrow \infty} \sim 8t^2 \exp(-2t^2)
$$
(in \cite{Kennedy}, Kennedy more precisely gives  the asymptotic expansion of this tail at all orders). Similar results are also available (see \cite{DW17}) for fragmentations related to the stable L\'evy trees. We refer to Example 4 in Section \ref{sec:exinfinite} for a discussion on those models.

%%%%%%%%%%%%%%%%%%%
\section{Main results}
\label{sec:mainresults}
%%%%%%%%%%%%%%%%%%%

%%%%%%%%%%%%%%%%%%%%
\subsection{On the tail of $\zeta$}
\label{sec:tail}
%%%%%%%%%%%%%%%%%%%%

As said, our objective is to give a non-logarithmic estimate of $\mathbb P(\zeta>t)$ as $t \rightarrow \infty$ in a fairly general framework.

\bigskip

\textbf{Main hypothesis.} Throughout the article we assume that
\begin{equation}
\label{hyp:main}
\tag{$\mathbf H$}
\limsup_{x \rightarrow \infty} \frac{\phi'(x)x}{\phi(x)}<1.
\end{equation}
This assumption is not very restrictive (note that for any $\phi$, $\phi'(x) \leq \phi(x)/x, \forall x \geq 0$ since $\phi$ is convex) and is verified by all `natural' examples of fragmentation processes that we may think off. In particular it holds as soon as $\phi$ is regularly varying at $\infty$ with some index  $\gamma \in [0,1)$.
(Note that when $\nu$ is finite, $\phi$ is regularly varying at $\infty$ with index $\gamma=0$.) In any case, note that since Hypothesis (\ref{hyp:main}) means that $\frac{\phi'(x)}{\phi(x)}\leq \frac{c}{x} $ for some $c \in (0,1)$ and all $x$ large enough, the two functions $\phi$ and $\psi$ are both bounded from above by multiple of power functions:
\begin{equation}
\label{majofunctions}
\phi(x) ~ \lesssim ~ x^c  \qquad \text{and} \qquad \psi(x) ~ \lesssim ~ x^{\frac{1}{1-c}}.
\end{equation}

\bigskip

Here is our main result:

\begin{theorem}
\label{thm:main}
Assume \emph{(\ref{hyp:main})}. Then
$$
\mathbb P(\zeta>t)  ~\asymp ~\left(\frac{\psi (|\alpha|t)}{t} \right)^{\frac{1}{|\alpha|}-1} \left( \psi'(|\alpha|t)\right)^{\frac{1}{2}} ~\exp \left(- \int_{x_{\psi}+1}^{t} \frac{\psi(|\alpha|r)}{\abs r} \mathrm dr \right).
$$ 
\end{theorem}

\bigskip

We note that the right-hand side only depends on the dislocation measure $\nu$ through the function $\phi$, which is known not to contain all the information about the dislocation of particles (two different dislocation measures may lead to the same $\phi$). We also emphasize that this right-hand side in fact only depends on the distribution under $\nu$ of the largest fragment (see the Remark in Section \ref{sec:intro}). To get a more explicit expression on specific examples, a basic idea is of course to find a sufficiently precise asymptotic expansion for $\phi$ and then deduce one for $\psi$. In all the examples we will develop in Section \ref{sec:examples} the above formula will reduced to the more pleasant form
$$
\mathbb P(\zeta>t)  ~\asymp ~ t^{a_0} \exp\left(-\sum_{i=1}^k a_i t^{\alpha_i}\right)
$$
for some explicit $a_i,\alpha_i$,
where the largest $\alpha_i$, say $\alpha_k$, is greater or equal to 1 (with equality iff $\nu$ is finite), with $a_k>0$ but the others $a_i$ of any sign. In all cases, the function in front of the exponential in the right-hand side of Theorem \ref{thm:main} is  bounded from above for large $t$ by power functions and the exponential term is smaller than $\exp(-at)$ for some $a>0$ (this is a consequence of (\ref{majofunctions}), and the fact that $\psi'(x) \lesssim \psi(x)/x$ under (\ref{hyp:main}) -- this is easy to see, see \cite[Lemma 6]{H20a} if needed -- and that $x \mapsto \psi(x)/x$ is increasing hence larger than a constant).

\medskip

Our strategy to prove Theorem \ref{thm:main} consists in comparing the behavior of the tail of $\zeta$ with that of the tail of the extinction time of a typical fragment  $\zeta_{\mathrm{tag}}$. This last time is strictly smaller than $\zeta$ and a key point is to understand how the ratio of their tails behaves asymptotically. We get the following:

\begin{prop}
\label{prop-eps-2}
Assume \emph{(\ref{hyp:main})}. Then,
$$
\mathbb P(\zeta >t)  ~\asymp ~\left(\frac{\psi (|\alpha|t)}{t} \right)^{\frac{1}{|\alpha|}} \mathbb P(\zeta_{\mathrm{tag}}>t).
$$
\end{prop}

\medskip

In particular we see that $$\mathbb P(\zeta >t)  \asymp  \mathbb P(\zeta_{\mathrm{tag}}>t) \quad \text{  iff } \quad \nu \text{ is finite}$$ since by (\ref{def:psi}) $\psi(\abs t) t^{-1} \rightarrow \abs \nu(\mathcal S^{\downarrow})$ as $t \rightarrow \infty$. Note that together with (\ref{eq:tailI}), Proposition \ref{prop-eps-2} immediately implies Theorem \ref{thm:main}. In fact, we expect that 
$$\mathbb P(\zeta >t)~ \propto ~ \left(\frac{\psi (\abs t)}{t} \right)^{\frac{1}{|\alpha|}} \mathbb P(\zeta_{\mathrm{tag}}>t),$$ 
which would immediately leads to an equivalence instead of bounds in Theorem \ref{thm:main}. However we do not now yet how to prove this equivalence, except in the following particular case:

\smallskip

\begin{cor}
\label{prop:nufinite}
Assume that $\int_{\mathcal S^{\downarrow}} (1-s_1)^{-1}\nu (\mathrm d \mathbf s)<\infty$ \emph{(}in particular $\nu$ is finite\emph{)}. Then, 
$$
\mathbb P(\zeta>t)  ~ \propto ~ \exp \big(-\nu(\mathcal S^{\downarrow}) t\big)
$$
with a constant of proportionality  $c:=\lim_{t \rightarrow \infty}\mathbb P(\zeta>t) \exp\big(\nu(\mathcal S^{\downarrow}) t\big)$ greater than 1.
\end{cor}

\bigskip

Last, coming back to a more general situation, we note that when $\phi$ is regularly varying at $\infty$ Theorem \ref{thm:main} leads to:

\begin{cor}
\label{prop:2}
Assume that $\phi$ is regularly varying at $\infty$ with index $\gamma \in [0,1)$. Then
$$
\mathbb P(\zeta>t)  ~\asymp ~\left(\frac{\psi (t)}{t} \right)^{\frac{1}{|\alpha|}-\frac{1}{2}} \exp \left(- \int_{x_{\psi}+1}^{t} \frac{\psi(\abs r)}{\abs r} \mathrm dr \right)
$$
where the function $x \mapsto \psi(x)/x$ is regularly varying at $\infty$ with index $\gamma/(1-\gamma)$.
\end{cor}

\bigskip

\textit{Remark.} If $\zeta_r$ designs the extinction time of an $(\alpha,r \nu)$-fragmentation, $r>0$, then $\zeta_r$ is distributed as $r^{-1} \zeta_1$. There will therefore be no loss of generality to consider a multiple of $\nu$ when necessary, when this simplifies the notation for example.

\bigskip

Proposition 3 will be proved in Section \ref{sec:ProofTh1}, while Corollary \ref{prop:nufinite} and Corollary \ref{prop:2} will be proved in Section \ref{sec:coro}.

%%%%%%%%%%%%%%%%%%%%%%%%%%%%%%%%%%%%%%%%%
\subsection{Asymptotics of moments of the largest and tagged fragments}
%%%%%%%%%%%%%%%%%%%%%%%%%%%%%%%%%%%%%%%%%

From Theorem \ref{thm:main} we deduce the large time behavior of all positive moments of the largest fragment $\mathbb E\left[F_1^a(t)\right]$, $a>0$. This is to be compared to the asymptotic behavior of moments of the tagged fragment, already known from previous works. Interestingly the moments of the largest fragment decrease significantly slower than those of the tagged fragment, but these two fragments conditioned on non extinction have an asymptotic behavior of the same order. 

\bigskip

\begin{prop}[Largest fragment]
\label{prop:decreaseF1}
Assume \emph{(\ref{hyp:main})}. Then for any $a>0$, 
\begin{eqnarray*}
\mathbb E\left[F_1^a(t)\right]  &~ \asymp ~& \left(\frac{t}{\psi (|\alpha|t)} \right)^{\frac{a}{|\alpha|}} \mathbb P(\zeta>t) \\
&~ \asymp ~& \left(\frac{t}{\psi (|\alpha|t)}\right)^{1+\frac{a-1}{|\alpha|}} \left( \psi'(|\alpha|t)\right)^{\frac{1}{2}} ~\exp \left(- \int_{x_{\psi}+1}^{t} \frac{\psi(|\alpha|r)}{\abs r} \mathrm dr \right).
\end{eqnarray*}
\end{prop}

\bigskip

To compare this with the behavior of the tagged fragment process $F_{\mathrm{tag}}$, we recall that this process is a non-increasing self-similar process. In \cite[Theorem 1.2]{HR12}  the large time behavior of such process conditioned on being positive was studied under the hypothesis (\ref{hyp:main}) on the Laplace exponent of the underlying subordinator  (see also  \cite[Theorem 3.1]{H10} for a version of this result under the more restrictive assumption that the Laplace exponent is regularly varying). In \cite[Corollary 3]{H20a} this was adapted to get the behavior of the positive moments:

\medskip

\begin{prop}[Tagged fragment - Corollary 3 of \cite{H20a}]
\label{prop:moments1}
Assume \emph{(\ref{hyp:main})}. Then for any $a>0$,
\begin{eqnarray*}
\mathbb E\big[F_{\mathrm{tag}}^a(t)\big] &~ \propto ~& \left(\frac{t}{\psi(\abs t)}\right)^{\frac{a}{\abs}}\mathbb P(\zeta_{\mathrm{tag}}>t) \\
& ~ \propto ~ &  \left(\frac{t}{\psi(\abs t)}\right)^{1+\frac{a}{\abs}} (\psi'(\abs t))^{1/2} \exp\left(- \int_{x_{\psi}+1}^{t} \frac{\psi(\abs r)}{\abs r} \mathrm dr\right).
\end{eqnarray*}
\end{prop}

\bigskip

We therefore observe that for all $a>0$
$$
\frac{\mathbb E\left[F_1^a(t)\right] }{\mathbb E\big[F_{\mathrm{tag}}^a(t)\big] } ~ \asymp ~ \frac{\mathbb P(\zeta>t)}{\mathbb P(\zeta_{\mathrm{tag}}>t)} ~ \asymp ~ \left(\frac{\psi(\abs t)}{t}\right)^{\frac{1}{\abs}}.
$$
The  behavior of $F_1(t)$ and $F_{\mathrm{tag}}(t)$ conditioned on non-extinction is thus of the same order, and more precisely for all $a>0$
$$
\mathbb E\big[F_1^a(t) | \zeta>t\big] ~ \propto ~   \left(\frac{t}{\psi(\abs t)}\right)^{\frac{a}{\abs}}  ~ \propto ~ \mathbb E\big[F_{\mathrm{tag}}^a(t) | \zeta_{\mathrm{tag}}>t\big].
$$

\bigskip

Note moreover that the behavior of the tagged fragment implies, together with Proposition \ref{prop-eps-2}, that for all $a \geq 1$,
$$
\mathbb E\Bigg[\sum_{i \geq 1} F^a_i(t) \Bigg] =\mathbb E\Big[F^{a-1}_{\mathrm{tag}}(t)\mathbf 1_{\{F_{\mathrm{tag}}(t)>0\}}\Big] ~ \asymp ~ \left(\frac{t}{\psi (|\alpha|t)} \right)^{\frac{a}{|\alpha|}} \mathbb P(\zeta>t).
$$ 

\vspace{0.7cm}

\textbf{Conjecture:} large-time behavior of $F$ conditioned on non-extinction. The above estimates, in part, lead us to conjecture that the following Yaglom limit should hold: 
\begin{equation*}
F(t) \left(\frac{\psi (|\alpha|t)}{t} \right)^{\frac{1}{|\alpha|}}  ~ | ~ \zeta>t \quad \text{ converges in distribution in $\ell_1^{\downarrow}$ to a non-trivial limit as }t \rightarrow \infty, 
\end{equation*}
where $\ell_1^{\downarrow}$ denotes the set of non-increasing, non-negative sequences with finite sum, endowed with the $\ell_1$ distance.
The limit should then be a quasi-stationary distribution for a version of $F$ extended to the space $\ell_1^{\downarrow}$. For comparison, we mention that Yaglom limits for non-monotonic growth self-similar fragmentation processes are established in \cite[Proposition 4.1 and Corollary 4.4]{BBCK18}.

%%%%%%%%%%%%%%%%%%%%%%%%%%%%%%%%%%%%%
\section{Proof of Proposition \ref{prop-eps-2}}
\label{sec:ProofTh1}
%%%%%%%%%%%%%%%%%%%%%%%%%%%%%%%%%%%%%

We recall that Proposition \ref{prop-eps-2} is a key ingredient to get Theorem \ref{thm:main} since this theorem is a direct consequence of Proposition \ref{prop-eps-2} and the estimate 
(\ref{eq:tailI}). In this section we prove Proposition \ref{prop-eps-2}. The general idea is to get upper and lower bounds for the tail $\mathbb P(\zeta>t)$ involving moments of one or two tagged fragments at time $t$ and then use the asymptotic behavior of these moments to conclude. In that aim, we start by introducing some notation related to two tagged fragments in Section \ref{sec:backtag}. Then we use in Section \ref{sec:compar} the first and second moment methods to compare  $\mathbb P(\zeta>t)$  with functions of moments of one or two tagged fragments at time $t$, and $t$. In Section \ref{sec:momentsjoint} we establish the asymptotic behavior of moments of two tagged fragments (for one tagged fragment this is already settled in Proposition \ref{prop:moments1}). Last in Section \ref{sec:conclusion} we put the pieces together to conclude.

%%%%%%%%%%%%%%%%%%%%%%%%%%%%
\subsection{Tagging two fragments}
\label{sec:backtag}
%%%%%%%%%%%%%%%%%%%%%%%%%%%%

We will need to tag independently two fragments and study their joint behavior. 

\bigskip

\textbf{Tagging two fragments.} Given the interval fragmentation $O_F$ associated to $F$ (see Section \ref{sec:intro}), consider two marked points $U_{(1)}$, $U_{(2)}$ both uniformly distributed on $(0,1)$, independent, and independent of $O_F$ and let $F_{\mathrm{tag},(k)}(t)$, $k \in \{1,2\}$, denote the size of the fragments containing $U_{(k)}$, $k \in \{1,2\}$, at time $t$ for all $t \geq 0$.  Of course $F_{\mathrm{tag},(1)}, F_{\mathrm{tag},(2)}$ are both distributed as $F_{\mathrm{tag}}$ in (\ref{lambda_sub}) and their extinction times are distributed as $\zeta_{\mathrm{tag}}$ (\ref{def_I}).
We will need for our study the asymptotic behavior as $t \rightarrow \infty$ of the (joint) moments $\mathbb E\big[F^a_{\mathrm{tag},(1)}(t)F^b_{\mathrm{tag},(2)}(t)\big]$ for $a,b > 0$. This is studied below in Section \ref{sec:momentsjoint}. 

\bigskip

\textbf{The separation time $\T$.} We let $\T$ be the first time at which the marked points $U_{(1)}$, $U_{(2)}$ belong to different fragments of $O_F$. Since $U_{(2)}$ is independent of $(O_F,U_{(1)})$ and uniformly distributed in $(0,1)$, we note that
\begin{equation}
\label{propT}
\mathbb P(\T>t ~ | ~F_{\mathrm{tag},(1)})= F_{\mathrm{tag},(1)}(t), \quad \forall t \geq 0,
\end{equation}
a key point in the following. Moreover, it is known (\cite{BertoinSSF}) that $O_F$ is Feller and so the fragmentation property holds for stopping times with respect to the natural filtration of this process. In particular it holds for the randomized stopping time $\T$, which implies the existence of two independent processes $\bar F_{\mathrm{tag},(1)}, \bar  F_{\mathrm{tag},(2)}$ both distributed as $F_{\mathrm{tag}}$ and independent of $(F_{\mathrm{tag},(1)},F_{\mathrm{tag},(2)},\T)$ such that for all $t \geq 0$
\begin{spreadlines}{12pt}
\begin{eqnarray}
\label{L1T}
&&\hspace{-0.72cm}\big(F_{\mathrm{tag},(1)}(t),F_{\mathrm{tag},(2)}(t)\big)\mathbf 1_{\{\T \leq t\}} ~ \overset{(\mathrm d)}= \\
\nonumber
&&\hspace{-1cm} \quad \big(F_{\mathrm{tag},(1)}(\T) \bar{F}_{\mathrm{tag},(1)}\big((t-\T)F^{\alpha}_{\mathrm{tag},(1)}(\T) \big), F_{\mathrm{tag},(2)}(\T) \bar{F}_{\mathrm{tag},(2)}\big((t-\T)F^{\alpha}_{\mathrm{tag},(2)}(\T) \big)\big)\mathbf 1_{\{\T \leq t\}}
\end{eqnarray}
\end{spreadlines}
(in fact there is an identity in distribution of \emph{processes}, i.e. simultaneously for all $t\geq 0$, but we will not need it). We will refer to this as the \emph{strong fragmentation property} applied at time $\T$. 

%%%%%%%%%%%%%%%%%%%%%%%%%%%%%%%%%%%%%%%%%%%%%%%%%%
\subsection{Connection between the tail of $\zeta$ and moments of tagged fragments}
\label{sec:compar}
%%%%%%%%%%%%%%%%%%%%%%%%%%%%%%%%%%%%%%%%%%%%%%%%%%

Using the first and second moments methods, we get:

\begin{prop}
\label{prop:encadrement}
Assume \emph{(\ref{hyp:main})}. Then for $t$ large enough
$$
\frac{\mathbb E\big[F_{\mathrm{tag}}(t)\big]^2}{\mathbb E\big[F_{\mathrm{tag},(1)}(t)F_{\mathrm{tag},(2)}(t)\big]} ~\leq ~ \mathbb P(\zeta>t) ~\lesssim  ~   \left(\frac{\psi(\abs t)}{t} \right)^{\frac{2}{\abs}} \mathbb E\big[F_{\mathrm{tag}}(t)\big] $$
and the right-hand side could be rewritten as
$$
\mathbb P(\zeta>t) ~\lesssim  ~    \left(\frac{\psi(\abs t)}{t} \right)^{\frac{1}{\abs}} \mathbb P(\zeta_{\mathrm{tag}}>t). $$
\end{prop}

\bigskip

\textbf{Proof.}
For $t\geq 0$, let
$$
S(t):= \sum_{i \geq 1} F^2_i(t).
$$
By definition of the tagged fragment processes, conditionally on $F$ and independently for $k \in \{1,2\}$,
$$F_{\mathrm{tag},(k)}(t)=F_i(t) \quad \text{with probability} \quad F_i(t), \quad i \geq 1,$$ 
and $F_{\mathrm{tag},(k)}(t)=0$ otherwise. Hence $\mathbb E[F_{\mathrm{tag},(k)}(t) | F]=S(t)$, $k \in \{1,2\}$ and $\mathbb E[F_{\mathrm{tag},(1)}(t) F_{\mathrm{tag},(2)}(t) | F]=(S(t))^2$, which gives
$$
\mathbb E\left[S(t)\right]=\mathbb E\big[F_{\mathrm{tag}}(t)\big] \quad \text{and} \quad \mathbb E\big[(S(t))^2 \big]=\mathbb E\big[F_{\mathrm{tag},(1)}(t)F_{\mathrm{tag},(2)}(t)\big].
$$

$\bullet$ This immediately leads to the lower bound of the proposition since $\{\zeta>t\}=\{S(t)>0\}$, which together with the second moment method gives
$$
\mathbb P(\zeta>t)=\mathbb P(S(t)>0) \geq \frac{\mathbb E[S(t)]^2}{\mathbb E[(S(t))^2]} = \frac{\mathbb E\big[F_{\mathrm{tag}}(t)\big]^2}{\mathbb E\big[F_{\mathrm{tag},(1)}(t)F_{\mathrm{tag},(2)}(t)\big]}.
$$

$\bullet$ To get the upper bound we use the identity (valid for all $t\geq 0$)
\begin{equation}
\label{eqzetaplus}
\left(\zeta-t\right)^+ ~\overset{(\mathrm d)}= ~\sup_{i \geq 1} F^{\abs}_i(t) \zeta^{(i)},
\end{equation}
where $(\zeta^{(i)})_{i \geq 1}$ is a sequence of i.i.d. random variables distributed as $\zeta$,  independent of $F(t)$. This is a simple consequence of the fragmentation property of $F$ at time $t$. Then, by Markov's inequality, we get for $t$ large enough
\begin{eqnarray*}
\mathbb P\left(\zeta>t+\frac{t}{\psi(\abs t)} \right) &=&\mathbb P\left(\sup_{i \geq 1}  F^{\abs}_i(t) \zeta^{(i)} >\frac{t}{\psi(\abs t)}\right) \\
&\leq& \left(\frac{\psi(\abs t)}{ t}\right)^{\frac{2}{\abs}} \mathbb E\left[\sum_{i\geq 1} F^2_i(t) \left(\zeta^{(i)}\right)^{\frac{2}{\abs}} \right]  \\
&=& \left(\frac{\psi(\abs t)}{ t}\right)^{\frac{2}{\abs}} \mathbb E\big[F_{\mathrm{tag}}(t)\big] \mathbb E\Big[\zeta^{\frac{2}{\abs}}\Big] \\
&\lesssim & \left(\frac{\psi(\abs t)}{ t}\right)^{\frac{1}{\abs}} \mathbb P(\zeta_{\mathrm{tag}}>t)
\end{eqnarray*}
where for the last line we use Proposition \ref{prop:moments1} and the fact that all positive moments of $\zeta$ are finite (since the tail of $\zeta$ decreases exponentially fast, at least). 
Now set for $t>0$ $$f(t):=t+\frac{t}{\psi(\abs t)} ~ \geq ~t$$ and note that
$\psi(\abs t)/t \leq \psi(\abs f(t))/f(t)$
since the function $t \mapsto \psi(\abs t)/t$ is increasing. Together with (\ref{extreme}), all this implies that for $t$ large enough
$$
\mathbb P\left(\zeta> f(t) \right) ~ \lesssim  ~ \left(\frac{\psi(\abs f(t))}{f(t)}\right)^{\frac{1}{\abs}}  ~ \mathbb P\left(\zeta_{\mathrm{tag}}>f(t) \right).
$$
To conclude and get the expected upper bound, it remains to note that the function $f$ is bijective from $(a,\infty)$ to $(f(a),\infty)$ for $a$ large enough. Indeed, clearly $f(t) \rightarrow \infty$ as $t \rightarrow \infty$ and
$$
f'(t)=1+\frac{1}{\psi(\abs t)}-\frac{t \abs\psi'(\abs t)}{(\psi(\abs t))^2}
$$
which tends to 1 as $t \rightarrow \infty$ since $\psi(t) \rightarrow \infty$ as $t \rightarrow \infty$ and the function $x \mapsto x \psi'(x)/\psi(x)$ is bounded for $x$ large enough (this is a consequence of (\ref{hyp:main}), already mentioned in Section \ref{sec:tail}). Finally, we indeed get that for $t$ large enough
$$
\mathbb P\left(\zeta> t \right) ~ \lesssim  ~ \left(\frac{\psi(\abs t)}{t}\right)^{\frac{1}{\abs}}  ~ \mathbb P\left(\zeta_{\mathrm{tag}}>t \right)
$$
or equivalently (by Proposition \ref{prop:moments1})
$$
 \mathbb P(\zeta>t) ~\lesssim  ~   \left(\frac{\psi(\abs t)}{t} \right)^{\frac{2}{\abs}} \mathbb E\big[F_{\mathrm{tag}}(t)\big].
 $$
$\hfill \square$

%%%%%%%%%%%%%%%%%%%%%%%%%%%%%%%%%%%%%
\subsection{Asymptotics of joint moments of two tagged fragments}
\label{sec:momentsjoint}
%%%%%%%%%%%%%%%%%%%%%%%%%%%%%%%%%%%%%

The aim of this section is to prove that for the two independently tagged fragments $F_{\mathrm{tag},(1)},F_{\mathrm{tag},(2)}$:

\begin{prop}
\label{prop:moments2}
For all $a,b>0$,
$$
\mathbb E\left[F^a_{\mathrm{tag},(1)}(t)F^b_{\mathrm{tag},(2)}(t)\right] ~\asymp ~\left(\frac{t}{\psi(\abs t)}\right)^{\frac{a+b+1}{\abs}}\mathbb P(\zeta_{\mathrm{tag}}>t).
$$
\end{prop}

\bigskip

The proof of the upper bound is quite technical and for this we will need the following intermediate result which gives an estimation by layers. We recall that $\T$ designs the separation time of the two tagged fragments.

\medskip

\begin{lemma}
\label{lem:decoup}
Fix $a,b>0$. For all $r \in (0,\infty)$, there exist $c_r, t_r \in (0,\infty)$ such that for all $~t \geq t_r$ and all integers $i \geq 1$
$$
\mathbb E\left[F^a_{\mathrm{tag},(1)}(t)F^b_{\mathrm{tag},(2)}(t)  \mathbf 1_{\left\{t-\frac{(i+1)t}{\psi(\abs t)} <T_{\mathrm{sep}} \leq t-\frac{it}{\psi(\abs t)}\right\}} \right] ~ \leq ~ \frac{c_r} {i^{r-\frac{1}{\abs}}} \left(\frac{t}{\psi(\abs t)} \right)^{\frac{a+b+1}{\abs}} \mathbb P(\zeta_{\mathrm{tag}}>t).
$$
\end{lemma}

\bigskip

With this lemma in hands and the conditional distribution (\ref{propT}), the proof of Proposition \ref{prop:moments2} holds as follows:

\bigskip

\textbf{Proof of Proposition \ref{prop:moments2}.} Fix $a,b>0$.

\smallskip

$\bullet$ To get the lower bound, we use (\ref{propT}) and that $F_{\mathrm{tag},(2)}(t)=F_{\mathrm{tag},(1)}(t)$ when $\T>t$:
\begin{eqnarray*}
\mathbb E\left[F^a_{\mathrm{tag},(1)}(t)F^b_{\mathrm{tag},(2)}(t)\right] &\geq& \mathbb E \left[F^a_{\mathrm{tag},(1)}(t)F^b_{\mathrm{tag},(2)}(t) \mathbf 1_{\{\T>t\}}\right] \\
&=& \mathbb E \left[F^{a+b}_{\mathrm{tag},(1)}(t)\mathbf 1_{\{\T>t\}}\right] \\
&=& \mathbb E \left[F^{a+b}_{\mathrm{tag},(1)}(t)\mathbb E \left[\mathbf 1_{\{\T>t\}}  | F_{\mathrm{tag},(1)}\right] \right] \\
&=& \mathbb E \left[F^{a+b+1}_{\mathrm{tag},(1)}(t) \right] \\
&\geq & c \left(\frac{t}{\psi(\abs t)}\right)^{\frac{a+b+1}{\abs}}\mathbb P(\zeta_{\mathrm{tag}}>t)
\end{eqnarray*}
for some $c \in (0,\infty)$ and all $t$ large enough, by Proposition \ref{prop:moments1}.

\bigskip

$\bullet$ To get the upper bound we write 
\begin{eqnarray*}
\mathbb E\left[F^a_{\mathrm{tag},(1)}(t)F^b_{\mathrm{tag},(2)}(t)\right] &=& \sum_{i=1}^{\infty}\mathbb E\left[F^a_{\mathrm{tag},(1)}(t)F^b_{\mathrm{tag},(2)}(t)  \mathbf 1_{\left\{t-\frac{(i+1)t}{\psi(\abs t)} <T_{\mathrm{sep}} \leq t-\frac{it}{\psi(\abs t)}\right\}} \right] \\
&+& \mathbb E\left[F^a_{\mathrm{tag},(1)}(t)F^b_{\mathrm{tag},(2)}(t)  \mathbf 1_{\left\{T_{\mathrm{sep}} > t-\frac{t}{\psi(\abs t)} \right\}}\right].
\end{eqnarray*}
By Lemma \ref{lem:decoup} (with e.g. $r=2+1/\abs$) the sum over $i \geq 1$ is smaller than a multiple of
$\left(\frac{t}{\psi(\abs t)} \right)^{\frac{a+b+1}{\abs}} \mathbb P(\zeta_{\mathrm{tag}}>t)$ 
for all $t$ large enough. Moreover, 
\begin{eqnarray*}
&& \mathbb E\left[F^a_{\mathrm{tag},(1)}(t)F^b_{\mathrm{tag},(2)}(t) \mathbf 1_{\left\{T_{\mathrm{sep}} > t-\frac{t}{\psi(\abs t)} \right\}}\right] \\
&\leq& \mathbb E\left[F^a_{\mathrm{tag},(1)}\left(t-\frac{t}{\psi(\abs t)}\right)F^b_{\mathrm{tag},(2)}\left(t-\frac{t}{\psi(\abs t)}\right) \mathbf 1_{\left\{T_{\mathrm{sep}} > t-\frac{t}{\psi(\abs t)} \right\}}\right] \\
&=& \mathbb E \left[F^{a+b+1}_{\mathrm{tag},(1)}\left(t-\frac{t}{\psi(\abs t)}\right) \right] \\
& \underset{\text{Prop. \ref{prop:moments1}}}\lesssim & \left(\frac{ t-\frac{t}{\psi(\abs t)}}{\psi\left( \abs \left( t-\frac{t}{\psi(\abs t)}\right)\right)} \right)^{\frac{a+b+1}{\abs}}\mathbb P\left(\zeta_{\mathrm{tag}}>t-\frac{t}{\psi(\abs t)}\right).
\end{eqnarray*}
We claim that this last upper bound is smaller than a multiple of $ \Big(\frac{t}{\psi(\abs t)} \Big)^{\frac{a+b+1}{\abs}} \mathbb P(\zeta_{\mathrm{tag}}>t)$ for $t$ large enough, as expected to conclude the proof of the lemma. To see this, recall that a consequence of (\ref{hyp:main}) is that $\psi'(x) \lesssim \psi(x)/x$ and note that this implies the existence of some $\kappa > 0$ (in fact $\kappa > 1$) such that $\psi(\lambda x) \leq \psi(x) \lambda^{\kappa}$ for all $\lambda \geq 1$ and all $x$ large enough (independent of $\lambda$). Hence
$$
\frac{ t-\frac{t}{\psi(\abs t)}}{\psi\left( \abs \left( t-\frac{t}{\psi(\abs t)}\right)\right)} ~ \lesssim ~ \frac{t}{\psi(\abs t)} .
$$
Besides, by (\ref{extreme}),
\begin{eqnarray*}
\mathbb P\left(\zeta_{\mathrm{tag}}>t-\frac{t}{\psi(\abs t)}\right) &\lesssim&  \mathbb P\left(\zeta_{\mathrm{tag}}>t-\frac{t}{\psi(\abs t)}+  \frac{2\left(t-\frac{t}{\psi(\abs t)}\right)}{\psi\left(\abs(t-\frac{t}{\psi(\abs t)})\right)}\right)  
\\
&\lesssim& \mathbb P\left(\zeta_{\mathrm{tag}}>t\right)
\end{eqnarray*}
since, recalling that $\psi$ is increasing and converges to $+\infty$,
\begin{eqnarray*}
t-\frac{t}{\psi(\abs t)}+  \frac{ 2\left(t-\frac{t}{\psi(\abs t)}\right)}{\psi\left(\abs(t-\frac{t}{\psi(\abs t)})\right)} &\geq& t -\frac{t}{\psi(\abs t)}+   \frac{2t}{ \psi(\abs t)}-  \frac{2t}{\psi(\abs t)\psi\left(\abs(t-\frac{t}{\psi(\abs t)}\right)} \\ 
&\geq & t \qquad \text{for $t$ large enough}.
\end{eqnarray*}
$\hfill \square$

\bigskip

It remains now to prove Lemma \ref{lem:decoup}. The rest of this section is devoted to this task. We start with a preliminary lemma in Section \ref{lem:prelim} and then turn to the proof of Lemma \ref{lem:decoup} in Section \ref{pr:lemdecoup}.

\subsubsection{A preliminary lemma}
\label{lem:prelim}

As noted in the previous section, the estimate
$$
\mathbb E\left[F^a_{\mathrm{tag},(1)}(t)\mathbf 1_{\{T_{\mathrm{sep}}>t\}} \right]  ~  \asymp ~ \left(\frac{t}{\psi(\abs t)}\right)^{\frac{a+1}{\abs}}\mathbb P(\zeta_{\mathrm{tag}}>t) \qquad \quad (a>0)
$$
is easy to settle. Here our goal is to complete this estimate by evaluating the expectation of $F^a_{\mathrm{tag},(1)}(t)$ on events where the separation time $T_{\mathrm{sep}}$ is slightly smaller than $t$. The correct scale to get useful information consists in considering the events $\left\{t\geq T_{\mathrm{sep}}>t- \frac{t(i+1)}{\psi(\abs t)} \right\}$, $i\geq 0$.

\begin{lemma}
\label{lem:prel1}
Fix $a>0$. Then  there exist  $c_0,t_0 \in (0,\infty)$ such that for all $t\geq t_0$ and all integers $i\geq 0$ 
$$
\mathbb E\left[F^a_{\mathrm{tag},(1)}(t)\mathbf 1_{\left\{T_{\mathrm{sep}}>t- \frac{t(i+1)}{\psi(\abs t)} \right\}} \right]
~\leq ~ c_0  (i+1)^{\frac{1}{\abs}}  \left(\frac{t}{\psi(\abs t)} \right)^{\frac{a+1}{\abs}} \mathbb P(\zeta_{\mathrm{tag}}>t). 
$$
\end{lemma}

\bigskip

\textbf{Proof.} Let $t_{\psi}\geq 1$ be some threshold large enough so that $\psi(\abs t)$ is well-defined for all $t \geq t_{\psi}$. 
Then fix $a>0$ and consider $t\geq t_{\psi}$. Note that if $i$ is an integer such that $t<t(i+1)/\psi(\abs t)$, then 
$$\mathbb E\left[F^a_{\mathrm{tag},(1)}(t)\mathbf 1_{\left\{T_{\mathrm{sep}}>t- \frac{t(i+1)}{\psi(\abs t)} \right\}} \right]=\mathbb E\left[F^a_{\mathrm{tag},(1)}(t)\right]$$
which is bounded from above by a multiple (independent of $i$) of $\left(\frac{t}{\psi(\abs t)} \right)^{\frac{a}{\abs}} \mathbb P(\zeta_{\mathrm{tag}}>t)$ by Proposition \ref{prop:moments1}, which it itself smaller than $$(i+1)^{\frac{1}{\abs}}  \left(\frac{t}{\psi(\abs t)} \right)^{\frac{a+1}{\abs}} \mathbb P(\zeta_{\mathrm{tag}}>t).$$ So in the following we will only consider integers $i$ such that $t \geq t(i+1)/\psi(\abs t)$. By (\ref{propT}) we can write
\begin{eqnarray*} 
\mathbb E\left[F^a_{\mathrm{tag},(1)}(t)\mathbf 1_{\left\{T_{\mathrm{sep}}>t- \frac{t(i+1)}{\psi(\abs t)} \right\}} \right]~ = ~\mathbb E \left[F^a_{\mathrm{tag},(1)}(t) F_{\mathrm{tag},(1)}\left(t-\frac{t(i+1)}{\psi(\abs t)}\right)\right].
\end{eqnarray*} 
To get the upper bound stated in the lemma, we split in two the expectation in the right-hand side of the above equality according to whether $$F^{\alpha}_{\mathrm{tag},(1)}\left(t-\frac{t(i+1)}{\psi(\abs t)}\right) \frac{t(i+1)}{\psi(\abs t)} \quad \text{is smaller or larger than 1}$$
and will bound from above each of the two expectations thus obtained. In that aim we use the fragmentation property at time $t-\frac{t(i+1)}{\psi(\abs t)}$ to write 
$$
F_{\mathrm{tag},(1)}(t) ~\overset{(\mathrm d)}=~F_{\mathrm{tag},(1)}\left(t-\frac{t(i+1)}{\psi(\abs t)}\right) \bar{F}_{\mathrm{tag},(1)}\left(F^{\alpha}_{\mathrm{tag},(1)}\left(t-\frac{t(i+1)}{\psi(\abs t)}\right) \frac{t(i+1)}{\psi(\abs t)} \right)
$$
where $\bar F_{\mathrm{tag},(1)}$ is a process independent of $F_{\mathrm{tag},(1)}$ and distributed as $F_{\mathrm{tag},(1)}$. This implies that for any $b>0$
\begin{eqnarray*}
\mathbb E\Big[F^{b}_{\mathrm{tag},(1)}(t)\Big] &\geq& \mathbb E\left[F^{b}_{\mathrm{tag},(1)}\left(t-\frac{t(i+1)}{\psi(\abs t)}\right)  \bar{F}^b_{\mathrm{tag},(1)}(1) \mathbf 1_{\left\{F^{\alpha}_{\mathrm{tag},(1)}\left(t-\frac{t(i+1)}{\psi(\abs t)}\right) \frac{t(i+1)}{\psi(\abs t)} \leq 1\right\}}\right] \\
&= & \mathbb E\Big[F_{\mathrm{tag},(1)}^b(1)\Big]  \mathbb E\left[F^{b}_{\mathrm{tag},(1)}\left(t-\frac{t(i+1)}{\psi(\abs t)} \right)  \mathbf 1_{\left\{F^{\alpha}_{\mathrm{tag},(1)}\left(t-\frac{t(i+1)}{\psi(\abs t)}\right) \frac{t(i+1)}{\psi(\abs t)} \leq 1\right\}}\right],
\end{eqnarray*}
and therefore, taking $b=a+1$ and using that $F_{\mathrm{tag},(1)}(t) \leq F_{\mathrm{tag},(1)}\Big(t-\frac{t(i+1)}{\psi(\abs t)}\Big)$
\begin{eqnarray*}
&& \mathbb E\left[F^a_{\mathrm{tag},(1)}(t) F_{\mathrm{tag},(1)}\left(t-\frac{t(i+1)}{\psi(\abs t)}\right) \mathbf 1_{\left\{F^{\alpha}_{\mathrm{tag},(1)}\left(t-\frac{t(i+1)}{\psi(\abs t)}\right) \frac{t(i+1)}{\psi(\abs t)} \leq 1\right\}}\right] \\
&\leq& \mathbb E\left[F_{\mathrm{tag},(1)}^{a+1}\left(t-\frac{t(i+1)}{\psi(\abs t)}\right) \mathbf 1_{\left\{F^{\alpha}_{\mathrm{tag},(1)}\left(t-\frac{t(i+1)}{\psi(\abs t)}\right) \frac{t(i+1)}{\psi(\abs t)} \leq 1\right\}}\right] \\
&\leq& \left(\mathbb E\Big[F_{\mathrm{tag},(1)}^{a+1}(1)\Big] \right)^{-1}\mathbb E\Big[F^{a+1}_{\mathrm{tag},(1)}(t)\Big] 
 \\
&\leq & d_1 \left(\frac{t}{\psi(\abs t)}\right)^{\frac{a+1}{\abs}}\mathbb P(\zeta_{\mathrm{tag}}>t) 
\end{eqnarray*}
 for some finite $d_1$ independent of $t$ large enough (and the integers $i$ such that $t \geq t(i+1)/\psi(\abs t)$), by Proposition \ref{prop:moments1}. On the other hand,  
\begin{eqnarray*}
 && \mathbb E\left[F^a_{\mathrm{tag},(1)}(t) F_{\mathrm{tag},(1)}\left(t-\frac{t(i+1)}{\psi(\abs t)}\right) \mathbf 1_{\left\{F^{\alpha}_{\mathrm{tag},(1)}\left(t-\frac{t(i+1)}{\psi(\abs t)}\right) \frac{t(i+1)}{\psi(\abs t)} > 1\right\}}\right] \\
 &\leq& \mathbb E\Big[F^a_{\mathrm{tag},(1)}(t)\Big] \left(\frac{t(i+1)}{\psi(\abs t)}\right)^{\frac{1}{\abs}} \\
 &\leq & d_2 \mathbb P(\zeta_{\mathrm{tag}}>t) \left(\frac{t}{\psi(\abs t)}\right)^{\frac{a+1}{\abs}} (i+1)^{\frac{1}{\abs}}
\end{eqnarray*}
for some finite $d_2$ independent of $t$ large enough (and the integers $i$ such that $t \geq t(i+1)/\psi(\abs t)$), again by Proposition \ref{prop:moments1}. Finally, both inequalities are available for all $t$ sufficiently large and we conclude by summing them. $\hfill \square$

\subsubsection{Proof of Lemma \ref{lem:decoup}}
\label{pr:lemdecoup}

The proof of Lemma \ref{lem:decoup} relies on the strong fragmentation property (\ref{L1T}), Lemma \ref{lem:prel1} and the (at-least) exponential decrease of $\mathbb P(\zeta_{\mathrm{tag}}>t)$. In the remaining, $a,b,r \in (0,\infty)$ are fixed.

\bigskip

As in the proof of the previous lemma we let $t_{\psi}\geq 1$ be some threshold large enough so that $\psi(\abs t)$ is well-defined for all $t \geq t_{\psi}$. Then, to lighten notations, we set for all $t\geq t_{\psi}$, $i \geq 1$,
$$
A_{t,i}:=\left\{t-\frac{(i+1)t}{\psi(\abs t)} <T_{\mathrm{sep}} \leq t-\frac{it}{\psi(\abs t)}\right\}.
$$
Note that $T_{\mathrm{sep}} \leq t$ on $A_{t,i}$. By the strong fragmentation property we have the existence of two independent processes $\bar  F_{\mathrm{tag},(1)}, \bar  F_{\mathrm{tag},(2)}$ both distributed as $F_{\mathrm{tag}}$ and independent of $(F_{\mathrm{tag},(1)},F_{\mathrm{tag},(2)},\T)$ such that
\begin{eqnarray*}
&& \mathbb E\left[F^a_{\mathrm{tag},(1)}(t) F^b_{\mathrm{tag},(2)}(t) \mathbf 1_{A_{t,i}}  \right] \\
&= & \mathbb E\left[ F^{a}_{\mathrm{tag},(1)}(T_{\mathrm{sep}})  \bar  F^a_{\mathrm{tag},(1)}\big((t-T_{\mathrm{sep}}) F^{\alpha}_{\mathrm{tag},(1)}(T_{\mathrm{sep}})\big)F^{b}_{\mathrm{tag},(2)}(T_{\mathrm{sep}}) \bar  F^b_{\mathrm{tag},(2)}\big((t-T_{\mathrm{sep}})F^{\alpha}_{\mathrm{tag},(2)}(T_{\mathrm{sep}})\big)\mathbf 1_{A_{t,i}}  \right] \\
&\leq & f_{(a,b)}(t,i) + f_{(b,a)}(t,i)
\end{eqnarray*}
where for $x,y>0$
$$
f_{(x,y)}(t,i) = \mathbb E\left[F^{x+y}_{\mathrm{tag},(1)}(T_{\mathrm{sep}})  \bar  F^x_{\mathrm{tag},(1)}\big((t-T_{\mathrm{sep}}) F^{\alpha}_{\mathrm{tag},(1)}(T_{\mathrm{sep}})\big) \bar F^y_{\mathrm{tag},(2)}\big((t-T_{\mathrm{sep}})F^{\alpha}_{\mathrm{tag},(1)}(T_{\mathrm{sep}})\big)  \mathbf 1_{A_{t,i}}  \right].
$$
The inequality is obtained by splitting the expectation according to whether $F_{\mathrm{tag},(1)}(T_{\mathrm{sep}}) \geq F_{\mathrm{tag},(2)}(T_{\mathrm{sep}})$ or not, and then using the exchangeability of $(F_{\mathrm{tag},(1)}(T_{\mathrm{sep}}),F_{\mathrm{tag},(2)}(T_{\mathrm{sep}}))$. Our goal now is to prove the existence of  a finite $c_r$ and a $t_r \geq t_{\psi}$ such that for all $t\geq t_r$ and all integers $i\geq 1$ 
\begin{equation}
\label{majo_f_a,b}
f_{(a,b)}(t,i) \leq  \frac{c_r} {i^{r-\frac{1}{\abs}}} \left(\frac{t}{\psi(\abs t)} \right)^{\frac{a+b+1}{\abs}} \mathbb P(\zeta_{\mathrm{tag}}>t)
\end{equation}
(which is sufficient to prove the lemma, by symmetry of the upper bound of (\ref{majo_f_a,b}) in $a,b$). 
In that aim we do another splitting, by writing
$$
A^{-}_{t,i}:=A_{t,i} \cap \{(t-T_{\mathrm{sep}})F^{\alpha}_{\mathrm{tag},(1)}(T_{\mathrm{sep}}) \leq t_{\psi}\} \quad \text{ and } \quad A^{+}_{t,i}:=A_{t,i} \cap \{(t-T_{\mathrm{sep}})F^{\alpha}_{\mathrm{tag},(1)}(T_{\mathrm{sep}}) > t_{\psi}\},
$$
$$
f^{-}_{(a,b)}(t,i) = \mathbb E\left[ F^{a+b}_{\mathrm{tag},(1)}(T_{\mathrm{sep}})  \bar  F^a_{\mathrm{tag},(1)}\big((t-T_{\mathrm{sep}}) F^{\alpha}_{\mathrm{tag},(1)}(T_{\mathrm{sep}})\big) \bar  F^b_{\mathrm{tag},(2)}\big((t-T_{\mathrm{sep}})F^{\alpha}_{\mathrm{tag},(1)}(T_{\mathrm{sep}})\big)  \mathbf 1_{A^{-}_{t,i}}  \right]
$$
and
$$
f^{+}_{(a,b)}(t,i) = \mathbb E\left[F^{a+b}_{\mathrm{tag},(1)}(T_{\mathrm{sep}})  \bar F^a_{\mathrm{tag},(1)}\big((t-T_{\mathrm{sep}}) F^{\alpha}_{\mathrm{tag},(1)}(T_{\mathrm{sep}})\big) \bar  F^b_{\mathrm{tag},(2)}\big((t-T_{\mathrm{sep}})F^{\alpha}_{\mathrm{tag},(1)}(T_{\mathrm{sep}})\big)  \mathbf 1_{A^{+}_{t,i}}  \right]
$$
so that $f_{(a,b)}(t,i)=f^{-}_{(a,b)}(t,i)+f^{+}_{(a,b)}(t,i)$. We will now set up suitable bounds for $f^{-}_{(a,b)}(t,i)$ and $f^{+}_{(a,b)}(t,i)$ in order to get (\ref{majo_f_a,b}).

\bigskip

$\bullet$ \textbf{Case 1: $\boldsymbol{(t-T_{\mathrm{sep}})F^{\alpha}_{\mathrm{tag},(1)}(T_{\mathrm{sep}})}$ is small}, which, roughly, implies that $T_{\mathrm{sep}}$ is not too small. 
On the event $A^{-}_{t,i}$ we have 
$$1 \leq t_{\psi}^r F^{|\alpha| r}_{\mathrm{tag},(1)}(T_{\mathrm{sep}})(t-T_{\mathrm{sep}})^{-r} \leq  t_{\psi}^r F^{|\alpha| r}_{\mathrm{tag},(1)}(T_{\mathrm{sep}}) \frac{\psi(\abs t)^r}{i^{r}t^{r}},$$ 
which leads us, together with the fact that $F^a_{\mathrm{tag},(1)}(u) \leq 1, F^b_{\mathrm{tag},(2)}(u) \leq 1$ for all $u \geq 0$, to:
\begin{equation*}
f^{-}_{(a,b)}(t,i) \leq   t_{\psi}^r  \frac{\psi(\abs t)^r}{i^r t^r}\mathbb E\left[F^{a+b+ \abs r}_{\mathrm{tag},(1)}(T_{\mathrm{sep}}) \mathbf 1_{A^{-}_{t,i}}  \right].
\end{equation*}
Besides, note that (still applying the strong fragmentation property)
\begin{eqnarray*}
\mathbb E\left[F^{a+b+ \abs r}_{\mathrm{tag},(1)}(t)\mathbf 1_{A_{t,i}}\right]
&=&  \left[F^{a+b+ \abs r}_{\mathrm{tag},(1)}(\T) \bar{F}^{a+b+ \abs r}_{\mathrm{tag},(1)}\big((t-\T)F^{\alpha}_{\mathrm{tag},(1)}(\T) \big)\mathbf 1_{A_{t,i}}\right] \\
&\geq & \left[F^{a+b+ \abs r}_{\mathrm{tag},(1)}(\T) \bar{F}^{a+b+ \abs r}_{\mathrm{tag},(1)}(t_{\psi})  \mathbf 1_{A_{t,i}^-}\right] \\
&=& \left[\bar{F}^{a+b+ \abs r}_{\mathrm{tag},(1)}(t_{\psi})\right]\left[F^{a+b+ \abs r}_{\mathrm{tag},(1)}(\T) \mathbf 1_{A_{t,i}^-} \right]
\end{eqnarray*}
and therefore 
$$
\mathbb E\left[ F^{a+b+ \abs r}_{\mathrm{tag},(1)}(T_{\mathrm{sep}}) \mathbf 1_{A^{-}_{t,i}}  \right] \leq d_r (i+1)^{\frac{1}{|\alpha|}} \left(\frac{t}{\psi(\abs t)}\right)^{\frac{a+b+\abs r+1}{\abs}} \mathbb P(\zeta_{\mathrm{tag}}>t) 
$$
by Lemma \ref{lem:prel1}, for some finite $d_r$ independent of $t$ large enough and $i\geq 1$. So finally we proved that for those $t,i$
\begin{equation}
\label{majofmoins}
 f^{-}_{(a,b)}(t,i) \leq d_r t_{\psi}  \frac{(i+1)^{\frac{1}{|\alpha|}}}{i^r}  \left(\frac{t}{\psi(\abs t)}\right)^{\frac{a+b+1}{\abs}}\mathbb P(\zeta_{\mathrm{tag}}>t).
\end{equation}

\bigskip

$\bullet$ \textbf{Case 2: $\boldsymbol{(t-T_{\mathrm{sep}})F^{\alpha}_{\mathrm{tag},(1)}(T_{\mathrm{sep}})}$  is large.} The proof is a bit more involved when we work on $A_{t,i}^+$. We will need the following observation: there exits a positive constant $c_{r,\mathrm{tag}}$ such that
\begin{equation}
\label{ref:majoF}
\bar G(x):=\mathbb P(\zeta_{\mathrm{tag}}>x) \leq \frac{c_{r,\mathrm{tag}}}{(\psi(\abs x))^r} \quad \text{for } x \geq t_{\psi}
\end{equation}
which results from the facts that $\bar G(x)$ decreases at least exponentially as $x \rightarrow \infty$ and that $\psi$ is bounded from above by a power function under the hypothesis (\ref{hyp:main}). Since ${(t-T_{\mathrm{sep}})F^{\alpha}_{\mathrm{tag},(1)}(T_{\mathrm{sep}})}>t_{\psi}$ on $A_{t,i}^+$ we can use this bound together with Proposition \ref{prop:moments1} to see that conditioning the expectation in the definition of $f^+_{(a,b)}(t,i)$ on $(T_{\mathrm{sep}}, F_{\mathrm{tag},(1)}(T_{\mathrm{sep}}))$ and then using the independence of $\bar F_{\mathrm{tag},(1)}, \bar F_{\mathrm{tag},(2)}$ leads to:
\begin{eqnarray*}
\label{ineg:FI}
 && f^+_{(a,b)}(t,i)  \\
 &\underset{\text{Prop. \ref{prop:moments1}}}\leq & d_1 \mathbb E\left[F^{a+b}_{\mathrm{tag},(1)}(T_{\mathrm{sep}})  \big(\bar G((t-T_{\mathrm{sep}})F^{\alpha}_{\mathrm{tag},(1)}(T_{\mathrm{sep}})) \big)^2 \left(\frac{F^{\alpha}_{\mathrm{tag},(1)}(T_{\mathrm{sep}})(t-T_{\mathrm{sep}})}{\psi\big(\abs F^{\alpha}_{\mathrm{tag},(1)}(T_{\mathrm{sep}})(t-T_{\mathrm{sep}}) \big)}\right)^{\frac{a+b}{\abs }} \mathbf 1_{A^{+}_{t,i}} \right] \\
\nonumber
&=& d_1 \mathbb E\left[ \big(\bar G((t-T_{\mathrm{sep}})F^{\alpha}_{\mathrm{tag},(1)}(T_{\mathrm{sep}})) \big)^2 \left(\frac{t-T_{\mathrm{sep}}}{\psi\big(\abs F^{\alpha}_{\mathrm{tag},(1)}(T_{\mathrm{sep}})(t-T_{\mathrm{sep}}) \big)}\right)^{\frac{a+b}{\abs }} \mathbf 1_{A^{+}_{t,i}} \right] \\
&\underset{(\ref{ref:majoF})} \leq & d_1  c_{r,\mathrm{tag}} \frac{\psi(\abs t)^r}{i^r t^r}\mathbb E\left[\bar G((t-T_{\mathrm{sep}})F^{\alpha}_{\mathrm{tag},(1)}(T_{\mathrm{sep}}))  \left(\frac{t-T_{\mathrm{sep}}}{\psi\left(\abs F^{\alpha}_{\mathrm{tag},(1)}(T_{\mathrm{sep}})(t-T_{\mathrm{sep}}) \right)}\right)^{\frac{a+b}{|\alpha|}+r} 
\mathbf 1_{A_{t,i}^+} \right]
\end{eqnarray*}
for some $d_1 \in (0,\infty)$ independent of $i$ and $t$, where for the second inequality we also used that $1 \leq (t-T_{\mathrm{sep}})^r\left(\frac{\psi(\abs t)}{it}\right)^r$ on the event $A_{t,i}^+~$ (in fact on the whole event $A_{t,i}$). To conclude, we write
\begin{eqnarray*}
 \mathbb E && \hspace{-1cm}\left[F^{a+b+\abs r}_{\mathrm{tag},(1)}(t) \mathbf 1_{A_{t,i}}\right]
= \mathbb E\left[F^{a+b+\abs r}_{\mathrm{tag},(1)}(\T) \bar{F}^{a+b+\abs r}_{\mathrm{tag},(1)}\big((t-\T) F^{\alpha}_{\mathrm{tag},(1)}(\T) \big) \mathbf 1_{A_{t,i}}\right] \\
 &\geq&   \mathbb E\left[F^{a+b+\abs r}_{\mathrm{tag},(1)}(\T)  \bar{F}^{a+b+\abs r}_{\mathrm{tag},(1)}\big((t-\T)F^{\alpha}_{\mathrm{tag},(1)}(\T) \big)\mathbf 1_{A^+_{t,i}} \right] \\
 &\geq&  d_2\mathbb E \left[F^{a+b+\abs r}_{\mathrm{tag},(1)}(\T) \bar G((t-\T)F^{\alpha}_{\mathrm{tag},(1)}(\T)) \left(\frac{(t-\T)F^{\alpha}_{\mathrm{tag},(1)}(\T)}{\psi(\abs (t-\T)F^{\alpha}_{\mathrm{tag},(1)}(\T))}\right)^{\frac{a+b}{\abs}+r} \mathbf 1_{A^+_{t,i}} \right] \\
&=&  d_2\mathbb E \left[\bar G((t-\T)F^{\alpha}_{\mathrm{tag},(1)}(\T)) \left(\frac{t-\T}{\psi(\abs (t-\T)F^{\alpha}_{\mathrm{tag},(1)}(\T))}\right)^{\frac{a+b}{\abs}+r} \mathbf 1_{A^+_{t,i}} \right]
\end{eqnarray*}
where the inequality between the second and third lines is obtained by conditioning on \linebreak $(\T, F_{\mathrm{tag},(1)}(\T))$ and using Proposition \ref{prop:moments1}: note that this proposition implies that  $\mathbb E\big[F^{a+b+\abs r}_{\mathrm{tag},(1)}(u)\big]$ is greater than $d_2 \left(\frac{t}{\psi(\abs u)}\right)^{\frac{a+b}{\abs}+r}\bar G(u)$ for some $d_2>0$ and all $u \geq t_{\psi}$ -- this constraint on $u$ not being too small is crucial here and that's partly why we restrict ourselves to $A^+_{t,i}$. Together with Lemma \ref{lem:prel1} which gives an upper bound for $\mathbb E\left[F^{a+b+\abs r}_{\mathrm{tag},(1)}(t)) \mathbf 1_{A_{t,i}}\right]$ this implies that
\begin{eqnarray*}
&& \mathbb E \left[\bar G((t-\T)F^{\alpha}_{\mathrm{tag},(1)}(\T)) \left(\frac{t-\T}{\psi(\abs (t-\T)F^{\alpha}_{\mathrm{tag},(1)}(\T))}\right)^{\frac{a+b}{\abs}+r} \mathbf 1_{A^+_{t,i}} \right] \\ 
&\leq& \frac{d_3}{d_2} (i+1)^{\frac{1}{\abs}} \left(\frac{t}{\psi(\abs t)}\right)^{\frac{a+b+1}{\abs}+r} \mathbb P(\zeta_{\mathrm{tag}}>t)
\end{eqnarray*}
for some $d_3$ independent of $i\geq 1$ and all $t$ large enough (the threshold being independent of $i$).
Finally we have proved that for $i \geq 1$ and all $t$ large enough
\begin{equation}
\label{majofplus}
 f^+_{(a,b)}(t,i) \leq \frac{d_1d_3 c_{r,\mathrm{tag}}}{d_2}   \frac{(i+1)^{\frac{1}{\abs}}}{i^r} \left(\frac{t}{\psi(\abs t)}\right)^{\frac{a+b+1}{\abs}} \mathbb P(\zeta_{\mathrm{tag}}>t).
\end{equation}

\bigskip

$\bullet$ \textbf{Conclusion.} We just add (\ref{majofmoins}) and (\ref{majofplus}) to get (\ref{majo_f_a,b}), which finishes the proof. $\hfill \square$

\bigskip

%%%%%%%%%%%%%%%%%%%%%%%%%%%%%%%%%%%%%
\subsection{Conclusion: end of the proof of Proposition \ref{prop-eps-2}}
\label{sec:conclusion}
%%%%%%%%%%%%%%%%%%%%%%%%%%%%%%%%%%%%%

The proof of Proposition \ref{prop-eps-2} is now easy to finish.

\bigskip

$\bullet$ From the right-hand inequality of Proposition \ref{prop:encadrement} and Proposition \ref{prop:moments1} we directly get that
$$
\mathbb P(\zeta > t)~= ~O\bigg(\left( \frac{\psi(\abs t)}{t}\right)^{\frac{1}{\abs}} \mathbb P(\zeta_{\mathrm{tag}}>t)\bigg). 
$$

\medskip

$\bullet$ From the left-hand inequality of Proposition \ref{prop:encadrement}, Proposition \ref{prop:moments1} and Proposition \ref{prop:moments2} we get 
\begin{eqnarray*}
\mathbb P(\zeta>t) &\geq&  \frac{\mathbb E\left[F_{\mathrm{tag}}(t)\right]^2}{\mathbb E\left[F_{\mathrm{tag},(1)}(t)F_{\mathrm{tag},(2)}(t)\right]} \\
&\gtrsim & ~ \frac{\left(\left(\frac{t}{\psi(\abs t)}\right)^{\frac{1}{\abs}}\mathbb P(\zeta_{\mathrm{tag}}>t)\right)^2}{\left(\frac{t}{\psi(\abs t)}\right)^{\frac{3}{\abs}}\mathbb P(\zeta_{\mathrm{tag}}>t).} \\
&=&  \left(\frac{\psi(\abs t)}{t}\right)^{\frac{1}{\abs}}\mathbb P(\zeta_{\mathrm{tag}}>t).
\end{eqnarray*}

%%%%%%%%%%%%%%%%%%%%%%%%%%%%%%%%%%%%%%%%%%
\section{Proof of Corollary \ref{prop:nufinite} and Corollary \ref{prop:2}}
\label{sec:coro}
%%%%%%%%%%%%%%%%%%%%%%%%%%%%%%%%%%%%%%%%%%

This short section is devoted to the proofs of these two corollaries stated in Section \ref{sec:tail}.

\bigskip

\textbf{Proof of Corollary \ref{prop:nufinite}.} We assume here that $\int_{\mathcal S^{\downarrow}} (1-s_1)^{-1}\nu (\mathrm d \mathbf s)<\infty$, in particular $\nu$ is finite. So we may and will also assume that $\nu(\mathcal S^{\downarrow})=1$, to ease notation. In this case, the fragmentation process $F$ remains constant equal to $(1,0,\ldots)$ during a time $T_1$ which has an exponential distribution with parameter 1 and then jumps to $F(T_1)$ which is distributed according to the probability $\nu$ and independent of $T_1$.  This, together with the fragmentation property at time $T_1$, leads to
\begin{eqnarray*}
\mathbb P(\zeta >t)&=& \mathbb P(T_1>t)+ \mathbb P\bigg(\sup_{i \geq 1} F_i^{\abs}(T_1) \zeta^{(i)}+T_1>t, T_1 \leq t \bigg) \\
&=& e^{-t} +\int_0^t e^{u-t}\mathbb P\bigg(\sup_{i \geq 1} F_i^{\abs}(T_1)\zeta^{(i)} \geq u\bigg) \mathrm du,
\end{eqnarray*}
where the $\zeta^{(i)},i\geq 1$ are i.i.d. distributed as $\zeta$ and independent of $(F(T_1),T_1)$. Hence $t \mapsto e^t \mathbb P(\zeta>t)$ is increasing on $\mathbb R_+$.
It remains to prove that it is also bounded on $\mathbb R_+$ (it will then be clear with the above equality that the limit is larger than 1).
Thanks to Proposition \ref{prop-eps-2} this is equivalent to show that $t \mapsto e^t \mathbb P(\zeta_{\mathrm{tag}}>t)$ is bounded on $\mathbb R_+$, since $\psi(\abs t)/t$ converges to $\abs \nu(\mathcal S^{\downarrow})$ here. And then, thanks to \cite[Section 5]{MZ06} (or \cite[Case 1 in Section 2.2]{H20a}), this is indeed the case as soon as the L\'evy measure $\pi$ of the subordinator $\abs \xi$ involved in the construction of $\zeta_{\mathrm{tag}}$ via (\ref{def_I}) satisfies $\int_0^{\infty} x^{-1} \pi(\mathrm dx)<\infty$. By (\ref{def:phi}), this L\'evy measure is defined in terms of $\nu$ by 
$$
\int_0^{\infty} f(x) \pi(\mathrm d x)=\int_{\mathcal S^{\downarrow}} \sum_{i \geq 1 }s_if(\abs |\ln(s_i)|) \nu(\mathrm d \mathbf s)
$$
for any positive measurable function $f$, and it is now easy to see that the finiteness of 
 $\int_0^{\infty} x^{-1} \pi(\mathrm dx)$ is equivalent to that of $\int_{\mathcal S^{\downarrow}} (1-s_1)^{-1}\nu (\mathrm d \mathbf s)$ (use that $|\ln(s_i)| \geq \ln (2)$ for all $i \geq 2$ and recall that $\int_{\mathcal S^{\downarrow}} (1-s_1)\nu (\mathrm d \mathbf s)<\infty$ by assumption).
$\hfill \square$

\bigskip

\textbf{Proof of Corollary \ref{prop:2}.} To deduce this corollary from Theorem \ref{thm:main} the only thing to check is that when $\phi$ is regularly varying the function $\psi$ is also regularly varying and $\psi'(x) \propto \psi(x)/x$.  This is a direct consequence of the definition of $\psi$ and on the fact that $\phi'$ is decreasing. See e.g. Lemma 7 in \cite{H20a} for details.
$\hfill \square$

%%%%%%%%%%%%%%%%%%%%%%%%%%%%%%%%%%%%%%%%%%%
\section{Proof of Proposition \ref{prop:decreaseF1}}
\label{sec:proofdecrease}
%%%%%%%%%%%%%%%%%%%%%%%%%%%%%%%%%%%%%%%%%%%

For the  proof of Proposition \ref{prop:decreaseF1} we will use Proposition \ref{prop-eps-2} and Proposition \ref{prop:moments1} several times. We proceed in two steps, each setting one of the inequalities. In the following we fix $a>0$.

\bigskip

$\bullet$ Let $p>1$ be such that $pa\abs>1$ and let $q:=p/(p-1)$. Using first H\"older's inequality and then the definition of the tagged fragment we obviously have that
\begin{eqnarray*}
\mathbb E\left[ F^{a\abs}_1(t)\right] = \mathbb E\left[ F^{a\abs}_1(t)\mathbf 1_{\{\zeta>t\}}\right] &\leq& \mathbb E\left[ F^{pa\abs}_1(t)\right]^{\frac{1}{p}} \mathbb P(\zeta>t)^{\frac{1}{q}} \\
&\leq &  \mathbb E\left[ \sum_{i=1}^{\infty} F^{pa\abs}_i(t)\right]^{\frac{1}{p}} \mathbb P(\zeta>t)^{\frac{1}{q}} \\
& = &  \mathbb E\left[ F^{pa\abs-1}_{\mathrm{tag}}(t)\right]^{\frac{1}{p}} \mathbb P(\zeta>t)^{\frac{1}{q}} \\
&\lesssim & \left( \frac{t}{\psi(\abs t)} \right)^{\frac{pa\abs-1}{\abs} \times \frac{1}{p}} \mathbb P(\zeta_{\mathrm{tag}}>t)^\frac{1}{p} \mathbb P(\zeta>t)^{\frac{1}{q}} 
\end{eqnarray*}
where the last inequality is obtained by Proposition \ref{prop:moments1}. Using Proposition \ref{prop-eps-2} and that $\frac{1}{p}+\frac{1}{q}=1$, the upper bound can be rewritten as the expected
$$
 \left( \frac{t}{\psi(\abs t)} \right)^{a}  \mathbb P(\zeta>t).
$$

\bigskip

$\bullet$ To get the inequality in the other direction, we will first show that
\begin{equation}
\label{boundzeta}
\mathbb E\Big[\big((\zeta-t)^+\big)^a\Big] ~\gtrsim ~ \left( \frac{t}{\psi(\abs t)} \right)^a \mathbb P(\zeta>t)
\end{equation}
where $x^+=\max(x,0)$. To prove this, we use a similar (but more precise) result for the random variable $\zeta_{\mathrm{tag}}$, which is that
\begin{equation}
\label{boundI}
\lim_{t \rightarrow \infty} \left(\frac{\psi(\abs t)}{t} \right)^a \frac{\mathbb E\left[\left((\zeta_{\mathrm{tag}}-t)^+\right)^a\right]}{\mathbb P(\zeta_{\mathrm{tag}}>t)} =\abs^{a} \Gamma(a+1).
\end{equation}
Roughly, this is a consequence of (\ref{extreme}). See  \cite[Section 6]{H20a} for a proof in the more general setting of exponential functionals of subordinators. Then use that
$$\mathbb E\Big[\big((\zeta-t)^+\big)^a\Big]=a \int_t^{\infty} \big((x-t)^+\big)^{a-1} \mathbb P(\zeta>x) \mathrm dx$$
together with Proposition \ref{prop-eps-2} to get that
$$
\mathbb E\Big[\big((\zeta-t)^+\big)^a\Big] ~ \asymp ~ a\int_t^{\infty} \left(\frac{\psi(\abs x)}{x}\right)^{\frac{1}{\abs}} \big((x-t)^+\big)^{a-1} \mathbb P(\zeta_{\mathrm{tag}}>x) \mathrm dx.
$$
Recalling that the function $x \mapsto \psi(x)/x$ is increasing, by definition of $\psi$ (see (\ref{def:psi})), the right-hand side of the last display can be bounded from below by
\begin{eqnarray*}
\left(\frac{\psi(\abs t)}{t}\right)^{\frac{1}{\abs}} a \int_t^{\infty} \big((x-t)^+\big)^{a-1} \mathbb P(\zeta_{\mathrm{tag}}>x) \mathrm dx &=& \left(\frac{\psi(\abs t)}{t}\right)^{\frac{1}{\abs}} \mathbb E\Big[\big((\zeta_{\mathrm{tag}}-t)^+\big)^a\Big] \\
& \underset{\text{by }(\ref{boundI})}\asymp & \left(\frac{\psi(\abs t)}{t}\right)^{\frac{1}{\abs}-a} \mathbb P(\zeta_{\mathrm{tag}}>t) \\
& \asymp & \left(\frac{t}{\psi(\abs t)}\right)^{a} \mathbb P(\zeta>t).
\end{eqnarray*}
Hence (\ref{boundzeta}).

\bigskip

Fix then $c \in \big(\frac{2}{1+a\abs},2\big)$ and recall the identity in distribution (\ref{eqzetaplus}) to get
\begin{eqnarray*}
\mathbb E\Big[\big((\zeta-t)^+\big)^{ca}\Big] 
& \leq & \mathbb E\left[ F^{\frac{ca\abs}{2}} _1(t) \sup_{i \geq 1} F^{\frac{ca\abs}{2}}_i(t) (\zeta^i)^{ca} \right] \\
& \leq &\mathbb E\left[ F^{a\abs}_1(t)\right]^{\frac{c}{2}}  \mathbb E\left[  \sup_{i \geq 1} F^{\frac{ca\abs}{2-c}}_i(t) (\zeta^i)^{\frac{2ca}{2-c}} \right]^{\frac{2-c}{2}}, 
\end{eqnarray*}
using again H\"older's inequality.  Next note that $\frac{ca\abs}{2-c}>1$ with our constraints on $c$ and then write
\begin{eqnarray*}
 \mathbb E\left[  \sup_{i \geq 1} F_i(t)^{\frac{ca\abs}{2-c}} (\zeta^i)^{\frac{2ca}{2-c}} \right] &\leq&  \mathbb E\left[  \sum_{i =1}^{\infty} F_i(t)^{\frac{ca\abs}{2-c}} (\zeta^i)^{\frac{2ca}{2-c}} \right]  \\
 &=&  \mathbb E\left[ F^{\frac{ca\abs}{2-c}-1}_{\mathrm{tag}}(t) \right]  \mathbb E\left[\zeta ^{\frac{2ca}{2-c}}\right] \\
 &\lesssim&   \left( \frac{t}{\psi(\abs t)} \right)^{\frac{ca}{2-c}-\frac{1}{\abs}} \mathbb P(\zeta_{\mathrm{tag}}>t) \\
 & \lesssim &  \left( \frac{t}{\psi(\abs t)} \right)^{\frac{ca}{2-c}} \mathbb P(\zeta>t) 
\end{eqnarray*}
the two last inequalities being consequences of Proposition \ref{prop:moments1} and Proposition \ref{prop-eps-2} respectively. 

\bigskip

Finally all this leads to
$$\left( \frac{t}{\psi(\abs t)} \right)^{ca} \mathbb P(\zeta>t) ~\lesssim ~  \mathbb E\Big[\big((\zeta-t)^+\big)^{ca}\Big] ~ \lesssim ~ \mathbb E\left[ F^{a\abs}_1(t)\right]^{\frac{c}{2}} \left( \frac{t}{\psi(\abs t)} \right)^{\frac{ca}{2}} \left(\mathbb P(\zeta>t) \right)^{\frac{2-c}{2}}
$$
and then 
$$
\left( \frac{t}{\psi(\abs t)} \right)^{a} \mathbb P(\zeta>t) ~ \lesssim ~ \mathbb E\left[ F^{a\abs}_1(t)\right],
$$
which finishes the proof.

%%%%%%%%%%%%%%%%%%%%%%%%
\section{Some examples and applications}
\label{sec:examples}
%%%%%%%%%%%%%%%%%%%%%%%%

Our objectif is to illustrate more concretely our main results  with simple or natural models of fragmentation processes, and related random trees. The general strategy to get explicit bounds is of course to find a sufficiently precise asymptotic expansion of the function $\phi$ (\ref{def:phi}) related to the model, then deduce one for $\psi$ (\ref{def:psi}) and conclude with Theorem \ref{thm:main} and/or its corollaries. To illustrate this we start with the following general result. A (simple) proof of the first part of it can be found in \cite[Lemma 19]{H20a}. The second part is a direct consequence of Corollary \ref{prop:2}. 

\begin{lemma}
\label{prop:aenuinfinite}
Assume that $\nu$ is a dislocation measure such that the corresponding $\phi$ writes
$$\phi(x)=x^{\gamma} \left(1-\sum_{i=1}^k \frac{c_i}{x^{\gamma_i}}+O(x^{-1-\varepsilon})\right)$$
for some $\gamma \in [0,1)$, $1/2<\gamma_1 <\gamma_2 \ldots <\gamma_{k-1} < \gamma_k = 1$, $\varepsilon>0$ and $c_i \in \mathbb R$, $1 \leq i \leq k$ . Then, 
$$
\frac{\psi(x)}{x}=x^{\frac{\gamma}{1-\gamma}}\left(1-\sum_{i=1}^k \frac{c_i}{(1-\gamma)x^{\frac{\gamma_i}{1-\gamma}}} \right)+O(x^{-1-\eta})
$$
for some $\eta>0$. Consequently for any $(\alpha,\nu)$-fragmentation, by Corollary \ref{prop:2},
\begin{eqnarray*}
\mathbb P(\zeta>t) &~\asymp ~& t^{\frac{\gamma}{1-\gamma}\left(\frac{1}{\abs}-\frac{1}{2}\right)+\frac{c_k}{\abs(1-\gamma)}}\exp\left(-\abs^{\frac{\gamma}{1-\gamma}} (1-\gamma)t^{\frac{1}{1-\gamma}}+\sum_{i=1}^{k-1} \frac{c_i \abs^{\frac{\gamma-\gamma_i}{1-\gamma}}}{1-\gamma_i} t^{\frac{1-\gamma_i}{1-\gamma}}\right) \\
&~\asymp ~& t^{\frac{\gamma}{1-\gamma} \frac{1}{\abs}} ~ \mathbb P(\zeta_{\mathrm{tag}} >t).
\end{eqnarray*}
\end{lemma}

\bigskip

From this, together with Proposition \ref{prop:decreaseF1} and Proposition \ref{prop:moments1}, we deduce a precise expression of the decrease of positive moments of the largest and tagged fragments. 
We emphasize that for the models where the function $\phi$ has an asymptotic expansion as above but with some $\gamma_i$s smaller or equal to $1/2$, there will be additional contributions to the non-$O(x^{-1-\eta})$ terms of $\psi(x)/x$, and consequently to the asympotic expression of $\mathbb P(\zeta>t)$: the closer a $\gamma_i$ is to 0, the more it contributes to a large number of terms (see for example the first calculations for the `Beta fragmentations' in the forthcoming Example 3).

\bigskip

Most of the examples that we detail below are based on the above lemma. We will also need the following asymptotic expansion for a quotient of Gamma functions: for all $c \in \mathbb R$,
\begin{equation}
\label{quo:Gamma}
\frac{\Gamma (x+c)}{\Gamma(x)}= x^{c}\left(1-\frac{c(1-c)}{2x}+O(x^{-2}) \right).
\end{equation}
This is a direct consequence of the asymptotic expansion of order 2 of the Gamma function $\Gamma(x)=e^{-x} x^{x-1/2} \sqrt{2\pi} \big( 1+(12x)^{-1}+O(x^{-2})\big).$ We will need as well the following lemma, which is a direct consequence of (\ref{quo:Gamma}).
\begin{lemma}
\label{lem:abBeta}
Let $a>0,b>-1$. Then
\begin{eqnarray*}
\frac{1}{\Gamma(b)}\int_0^1(1-u^x)u^{a-1}(1-u)^{b-1} \mathrm du &=& \frac{\Gamma(a)}{\Gamma(a+b)}- \frac{\Gamma(x+a)}{\Gamma(x+a+b)} \\
&\underset{\text{as $x \rightarrow \infty$}}{=}& \frac{\Gamma(a)}{\Gamma(a+b)}-  x^{-b}+b\left(a+\frac{b}{2}-\frac{1}{2}\right)x^{-b-1}+O(x^{-b-2})
\end{eqnarray*}
where $1/\Gamma$ denotes the extension to $\mathbb C$ by analytic continuation of the function $1/\Gamma$ initially defined on $\{z\in \mathbb C:\mathrm{Re}(z)>0\}$.
\end{lemma}

\subsection{Examples with finite splitting rate}

When the dislocation measure $\nu$ is finite, the quotient $\psi(\abs t)/t$ converges to a non-zero limit as $t \rightarrow \infty$ and the tails $\mathbb P(\zeta>t)$ and $\mathbb P(\zeta_{\mathrm{tag}} >t)$ are therefore asymptotically proportional (by Proposition \ref{prop-eps-2}) and also proportional to the moments of the largest and tagged fragments (by Proposition \ref{prop:decreaseF1}  and  Proposition \ref{prop:moments1} respectively). We can thus restrict ourselves here to statements on the tails of the extinction times. 

\medskip

Since there is no loss of generality in considering a multiple of the measure $\nu$, we will assume throughout this section that  $\nu$ is a probability. Our goal is to show with simple examples how the distribution of the largest piece under $\nu$ influences the asymptotic behavior of $\mathbb P(\zeta>t)$. We start with the following immediate consequence of Corollary \ref{prop:nufinite}.

\bigskip

\textbf{Ex.1. Fragmentations into $k$ identical pieces.} In these cases a fragment of mass $m$ splits into $k$ fragments of same mass $m/k$. Then $\nu=\delta_{(\frac{1}{k},\ldots,\frac{1}{k})}$, and for any $\alpha<0$ the extinction time $\zeta$ of an $(\alpha,\nu)$-fragmentation satisfies
$$
\mathbb P(\zeta>t) ~ \propto ~ \mathrm  \exp(-t).$$

\bigskip

\textbf{Ex.2. Uniform fragmentation into $k$ pieces.}  Here a fragment splits uniformly into $k$ sub-fragments, i.e. according to a Dirichlet $\mathrm{Dir}(1,\ldots,1)$ distribution on the $k-1$-dimensional simplex. The measure $\nu$ is then the push-forward of this distribution by the decreasing rearrangement function and 
$\phi(x)=1- k \mathbb E\big[B_k^{x+1}\big]$ where $B_k$ follows a $\mathrm{Beta}(1,k-1)$-distribution. Hence 
$$
\phi(x)=1-\frac{\Gamma(k+1)\Gamma(x+2)}{\Gamma(k+x+1)}, \quad \forall x \geq 0.
$$
The behavior of the tail of the extinction time of such an $(\alpha,\nu)$-fragmentation differs therefore significantly according to whether $k=2$ or $k \geq 3$. In the latter case, $\phi(x)=1+O(x^{-2})$ whereas when $k=2$, $\phi(x)=x/(x+2)$ (and therefore $\psi(x)=x-2$ for $x \geq 2$). This observation, together with Lemma \ref{prop:aenuinfinite} (or directly Corollary \ref{prop:2}) leads to 
$$
\mathbb P(\zeta>t) ~\asymp~t^{\frac{2}{\abs}} \exp(- t) \quad \text{when }k=2, \qquad \mathbb P(\zeta>t) ~\asymp~ \exp(- t) \text{ when }k \geq 3.
$$

\bigskip

\textbf{Ex.3.  Beta fragmentations.}  To understand better the role of the distribution of the largest piece under the probability $\nu$, consider more generally  the (binary) cases where a fragment of mass $m$ splits into two fragments of masses $mB$ and $m(1-B)$ where $B$ has a $\mathrm{Beta}(a,b)$ distribution on $[0,1]$, $a,b>0$. By symmetry we can assume that $b\geq a$. (Note that when $a=b=1$ we recover the uniform binary splitting.)   The measure $\nu$ is then the distribution of $(\max(B,1-B),\min(B,1-B))$ and
$$
\phi(x)=1-\mathbb E[B^{x+1}]-\mathbb E[(1-B)^{x+1}]=1- \frac{\Gamma(a+b)\Gamma(x+a+1)}{\Gamma(a)\Gamma(x+a+b+1)}-\frac{\Gamma(a+b)\Gamma(b+x+1)}{\Gamma(b)\Gamma(x+a+b+1)}
$$ 
so by (\ref{quo:Gamma}),
$$\phi(x)=1-\frac{\Gamma(a+b)}{\Gamma(b)x^{a}}-\frac{\Gamma(a+b)}{\Gamma(a)x^{b}}+O(x^{-a-1}).$$
When $a>1/2$, Lemma \ref{prop:aenuinfinite} gives the asymptotic behavior of $\mathbb P(\zeta>t)$ for such an $(\alpha,\nu)$-fragmentation, leading to six different situations:
\renewcommand{\arraystretch}{1.5}
$$
\mathbb P(\zeta>t) ~ \asymp ~\exp(-t) \times \left\{\begin{array}{ll} 1&   \text{ if } b\geq a >1  \\ t^{\frac{b}{\abs}}&   \text{ if } b> a =1 \\ t^{\frac{2}{\abs}} & \text{ if } b=a =1 \\ \exp\left(\frac{\Gamma(a+b)t^{1-a}}{\Gamma(b)(1-a)\abs ^a}\right) &  \text{ if } b> 1>a>1/2   \\\exp\left(\frac{\Gamma(a+b)t^{1-a}}{\Gamma(b)(1-a)\abs^a}\right) t^{\frac{a}{\abs}} &  \text{ if } 1=b > a>1/2 \\ \exp\left(\frac{\Gamma(a+b)t^{1-a}}{\Gamma(b)(1-a)\abs^a}+\frac{\Gamma(a+b)t^{1-b}}{\Gamma(a)(1-b)\abs^b}\right) &\text{ if } 1>b \geq a>1/2.  \end{array}\right.
$$
Otherwise, when $a \leq 1/2$, the smaller $a$ and $b$ are, the more terms there are in the estimates of $\mathbb P(\zeta>t)$. To illustrate, we consider the cases when $b \in (1/2,1)$ and $a \in (1/3,1/2]$. It is easy to see that then:
$$
\frac{\psi(x)}{x}=1-\frac{\Gamma(a+b)}{\Gamma(b)x^a}-\frac{\Gamma(a+b)}{\Gamma(a)x^b}-\frac{a(\Gamma(a+b))^2}{(\Gamma(b))^2x^{2a}}-\frac{(a+b)(\Gamma(a+b))^2}{\Gamma(a)\Gamma(b)x^{a+b}}+O(x^{-2b})+O(x^{-3a})+O(x^{-a-1}).
$$
And so for $b \in (1/2,1)$:
\renewcommand{\arraystretch}{1.5}
\begin{eqnarray*}
\mathbb P(\zeta>t) &~ \asymp ~& \exp\left(-t+\frac{\Gamma(a+b)t^{1-a}}{\Gamma(b)(1-a)\abs^a}+\frac{\Gamma(a+b)t^{1-b}}{\Gamma(a)(1-b)\abs^b}\right) \\  
&&\times \left\{\begin{array}{ll}  t^{\frac{(\Gamma(b+1/2))^2}{2(\Gamma(b))^2\abs}} &   \text{ if } a=1/2 \\ \exp\left(\frac{a(\Gamma(a+b))^2t^{1-2a}}{(\Gamma(b))^2(1-2a)\abs^{2a}}\right)& \text{ if }  a \in (1/3,1/2), ~a+b >1 \\ \exp\left(\frac{a(\Gamma(a+b))^2t^{1-2a}}{(\Gamma(b))^2(1-2a)\abs^{2a}}\right) t^{\frac{1}{\Gamma(a)\Gamma(b)\abs}}&   \text{ if }  a \in (1/3,1/2), ~a+b=1 \\  \exp\left(\frac{a(\Gamma(a+b))^2t^{1-2a}}{(\Gamma(b))^2(1-2a)\abs^{2a}}+\frac{(a+b)(\Gamma(a+b))^2t^{1-a-b}}{\Gamma(a)\Gamma(b)(1-a-b)\abs^{a+b}}\right)&   \text{ if } a \in (1/3,1/2), ~a+b <1.   \end{array}\right.
\end{eqnarray*}

\subsection{Examples with infinite splitting rate} 
\label{sec:exinfinite}

When $\nu$ is infinite we illustrate our results with self-similar fragmentations related to three families of random rooted real trees: (i) the stable L\'evy trees of Duquesne, Le Gall and Le Jan, (ii) the scaling limits of the so-called Ford's alpha model and (iii) of the Aldous beta-splitting model. (These families all contain as a fundamental example the Brownian tree of Aldous.) It turns out that these cases are mostly in the framework of Lemma \ref{prop:aenuinfinite}. Other examples of random trees related to self-similar fragmentations, as e.g. the $k$-ary trees studied in \cite{HS14} can be treated with similar calculations. We leave the details to the interested reader.

\bigskip

In fact, any fragmentation process with a negative index of self-similarity is associated  to a random compact rooted real tree,  equipped with a probability measure, that describes its genealogies \cite{HM04}. The extinction time $\zeta$ of the process can then be interpreted as the height (the maximal distance to the root) of the associated real tree and the extinction time of a tagged fragment, $\zeta_{\mathrm{tag}}$, as the height of a typical leaf (chosen according to the probability measure on the tree). As studied in \cite{HMPW08, HM12}, these self-similar fragmentation trees are the scaling limits (in distribution, with respect to the Gromov-Hausdorff-Prokhorov topology) of sequences of discrete trees, satisfying a so-called \emph{Markov-Branching property}. In particular the rescaled heights of these discrete trees converge to the extinction time  of the corresponding self-similar fragmentation, for which we have now precise information on the tail.

\bigskip

\textbf{Ex.4. Stable fragmentations and height of stable L\'evy trees.} For $\gamma \in (1,2]$, we let $\mathcal T_{\gamma}$ denote a stable L\'evy tree with branching mechanism $x \mapsto x^{\gamma}$. This tree is a compact rooted real tree that describes the genealogy of a continuous-state branching process with branching mechanism $x \mapsto x^{\gamma}$ \cite{LGLJ98, DLG02} and arises as the scaling limit of Galton-Watson trees conditioned on having $n$ vertices and with an offspring distribution in the domain of attraction of a $\gamma$-stable distribution \cite{Duq03}. It is naturally equipped with a probability measure $\mu_{\gamma}$ which is the limit of the uniform measure on the vertices of the conditioned Galton-Watson trees. In particular, the tree $\mathcal T_2$ is a multiple of the Brownian tree of Aldous. 

\bigskip

Bertoin \cite{BertoinSSF} for the Brownian tree and then Miermont \cite{M03} for the whole family of stable L\'evy trees shew that one gets a self-similar fragmentation process by considering the $\mu_{\gamma}$-masses of the connected components of the set of points of $\mathcal T_{\gamma}$ having a height greater than $t$, $t \geq 0$.  Its index of self-similarity is $\gamma^{-1}-1$ and the function $\phi$ associated to its dislocation measure is given by
$$
\phi(x)=\gamma \frac{\Gamma (x+1-\frac{1}{\gamma})}{\Gamma(x)}=\gamma x^{1-\frac{1}{\gamma}}\left(1-\frac{(1-\frac{1}{\gamma})\frac{1}{\gamma}}{2x}+O(x^{-2}) \right).
$$

Below we denote by $F_{\gamma,1}(t)$ the mass of the largest fragment of such a fragmentation at time $t$, by $F_{\gamma,\mathrm{tag}}(t)$ the mass of a tagged fragment, by  
$\zeta_{\gamma}$ the height of $\mathcal T_{\gamma}$ and by $\zeta_{\gamma,\mathrm{tag}}$ the height of a typical leaf of $\mathcal T_{\gamma}$, $\gamma \in (1,2]$. With our notation of Section \ref{sec:intro}, $\zeta_2$ is then distributed as $\sqrt 2 ~\zeta_{\mathrm{Br}}$.

\bigskip

Then, by Lemma \ref{prop:aenuinfinite}, we obtain:
\begin{cor} For all $\gamma \in (1,2]$,
$$\mathbb P(\zeta_{\gamma}>t) ~ \asymp ~  t^{1+\frac{\gamma}{2}}\exp\big(- (\gamma-1)^{\gamma-1} t^{\gamma}\big) \quad \text{and} \quad \mathbb P(\zeta_{\gamma, \mathrm{frag}}>t) ~ \asymp ~ t^{1-\frac{\gamma}{2}}\exp\big(- (\gamma-1)^{\gamma-1} t^{\gamma}\big).$$ 
Moreover for all $a>0$
$$
\hspace{-0.5cm}\mathbb E\left[F^{a}_{\gamma,1}(t) \right] ~ \asymp ~ t^{1+\gamma (\frac{1}{2}-a)} \exp\big(- (\gamma-1)^{\gamma-1} t^{\gamma}\big) ~ \asymp ~ t^{\gamma} \mathbb E\left[F^a_{\gamma,\mathrm{tag}}(t) \right].
$$
\end{cor}

This estimate of the tail of the height $\zeta_{\gamma}$ of the three $\mathcal T_{\gamma}$ is not new: in \cite[Theorem 1.5]{DW17}, Duquesne and Wang obtain a much more precise result by giving the  asymptotic expansion of this tail at all orders. We just wanted to emphasize that we easily get the first order as a consequence of our approach, and moreover some precise information on the decrease of the moments of the largest fragment. 

\bigskip

\textbf{Ex.5. Height of scaling limits of Ford's trees.} Ford \cite{Ford05} introduced a family of random cladograms (i.e. rooted binary finite trees) as possible models for phylogenetic trees. It is a family, indexed by a parameter $\alphaF \in [0,1]$, of sequences of growing trees constructed algorithmically according to a simple rule that only depends on the parameter $\alphaF$, that interpolates between the Yule tree ($\alphaF=0$), the uniform binary tree ($\alphaF=1/2$) and the comb tree ($\alphaF=1$). See e.g. \cite{CFW09, CDKM15, NW20, St09} for a few probabilistic references related to these models. For $\alphaF \in (0,1)$,  as the number of leaves of the trees of the sequence indexed by $\alphaF$ tends to infinity, one observes in the scaling limit a self-similar fragmentation tree (\cite[Section 5.2]{HMPW08}) with index of self-similarity $-\alphaF$ and a \emph{binary} dislocation measure (i.e. $\nu(s_1+s_2<1)=0$) characterized by
$$
\nu (s_1 \in \mathrm d u) =\frac{1}{\Gamma (1-\alphaF)}\left(\alphaF (u(1-u))^{-\alphaF-1}+(2-4\alphaF) (u(1-u))^{-\alphaF}\right), \quad u \in (1/2,1).
$$
Let $\zeta_{\alphaF}$ denote the height of this tree, $\zeta_{\alphaF,\mathrm{tag}}$ the height of a typical leaf and $F_{\alphaF}$ the corresponding fragmentation process.
(We emphasize that when $\alphaF=1/2$ this tree corresponds to a multiple of the Brownian tree.) The above definition of $\nu$ leads to
\begin{eqnarray*}
\phi(x)
= \frac{\alphaF}{\Gamma (1-\alphaF)} \int_0^1 (1-u^{x}) u^{-\alphaF}(1-u)^{-\alphaF-1}  \mathrm du+\frac{2-4\alphaF}{\Gamma (1-\alphaF)}  \int_0^1 (1-u^{x}) u^{1-\alphaF}(1-u)^{-\alphaF}  \mathrm du
\end{eqnarray*}
and then, by Lemma \ref{lem:abBeta}, to
\begin{eqnarray*}
\phi(x)
&=& x^\alphaF \left(1 - \left(\frac{3}{2}\alphaF^2-\frac{9}{2}\alphaF+2\right)x^{-1}+O(x^{-2}) \right).
\end{eqnarray*}
By Lemma \ref{prop:aenuinfinite}, we can therefore conclude that:
\begin{cor} For $\texttt{a} \in ]0,1[$,
\begin{eqnarray*}
\mathbb P(\zeta_{\texttt{a}}>t) 
~ \asymp ~ t^{\frac{2\texttt{a}^2-7\texttt{a}+4}{2\texttt{a}(1-\texttt{a})}} \exp\big(-\texttt{a}^{\frac{\texttt{a}}{1-\texttt{a}}}(1-\texttt{a})t^{\frac{1}{1-\texttt{a}}}\big)  ~ \asymp ~  t^{\frac{1} {1-\texttt{a}}} \cdot \mathbb P(\zeta_{\texttt{a},\mathrm{tag}}>t).
\end{eqnarray*}
Moreover for all $r>0$,
\begin{equation*}
\mathbb E\big[\big(F_{\texttt{a},1}(t)\big)^r \big] ~ \asymp ~ t^{-\frac{r}{1-\texttt{a}}} \cdot  t^{\frac{2\texttt{a}^2-7\texttt{a}+4}{2\texttt{a}(1-\texttt{a})}} \exp\big(-\texttt{a}^{\frac{\texttt{a}}{1-\texttt{a}}}(1-\texttt{a})t^{\frac{1}{1-\texttt{a}}}\big) ~ \asymp ~  t^{\frac{1} {1-\texttt{a}}} \cdot \mathbb E\big[\big(F_{\texttt{a},\mathrm{tag}} (t)\big)^r\big].
\end{equation*}
\end{cor}

\bigskip

\textbf{Ex.6. Height of $\beta$-splitting trees.} By $\beta$-splitting trees we refer to a family of random rooted real trees $(\mathcal T_{\beta})$ indexed by a parameter $\beta \in (-2,-1)$ that arise as the scaling limits (see \cite[Section 5.1]{HMPW08}) of well-known models of cladograms  introduced by Aldous in \cite{Ald96}. These real trees belong to the family of self-similar fragmentation trees and naturally extend the Beta-fragmentations consider earlier by replacing the Beta-dislocation `probability' measure by an `infinite $\sigma$-finite' Beta-dislocation measure with parameters $(\beta+1,\beta+1)$. More precisely, $\mathcal T_{\beta}$ has for index of self-similarity $\alpha=1+\beta$ and a \emph{binary} dislocation measure defined by
$$
\nu (s_1 \in \mathrm d u) =\frac{-1}{\Gamma (1+\beta)}(u(1-u))^{\beta}, \quad u \in (1/2,1).
$$
The tree $\mathcal T_{-3/2}$ is a multiple of the Brownian tree.
By Lemma \ref{lem:abBeta} we notice that
\begin{eqnarray*}
\phi(x)&=&\frac{-1}{\Gamma (1+\beta)} \int_0^1 (1-u^x) u^{\beta+1}(1-u)^{\beta}\mathrm du \\
&=& x^{-\beta-1}\left(1-\frac{\Gamma(\beta+2)}{\Gamma(2\beta+3)}x^{\beta+1}-\frac{(\beta+1)(3\beta+4)}{2}x^{-1} +O(x^{-2})\right).
\end{eqnarray*}
Then by Lemma \ref{prop:aenuinfinite} we get ($\zeta_{\beta}$ denotes the height of the tree, $\zeta_{\beta, \mathrm{tag}}$ the height of a typical leaf, and $F_{\beta}$ the associated fragmentation process):
\begin{cor}
For $\beta \in (-2,-3/2]$
$$
\mathbb P(\zeta_{\beta}>t)  ~ \asymp ~ t^{\frac{-2\beta-1}{2(\beta+2)}} \exp\Big(-a_{\beta} t^{\frac{1}{\beta+2}}+b_{\beta} t\Big) ~ \asymp ~ t^{\frac{1}{\beta+2}} \cdot \mathbb P(\zeta_{\beta,\mathrm{tag}}>t)  
$$
where
$$
a_{\beta}=(-\beta-1)^{\frac{-\beta-1}{\beta+2}}(\beta+2) \qquad \text{and} \qquad b_{\beta}=\frac{(2\beta+3) \Gamma(\beta+2)}{(\beta+2) \Gamma(2\beta+4)}.
$$
Moreover for all $a>0$
$$
\mathbb E\left[F^{a}_{\beta,1}(t) \right] ~ \asymp ~ t^{-\frac{a}{2+\beta}}  t^{\frac{-2\beta-1}{2(\beta+2)}} \exp\Big(-a_{\beta} t^{\frac{1}{\beta+2}}+b_{\beta} t\Big) ~ \asymp ~ t^{\frac{1}{2+\beta}} \mathbb E\left[F^{a}_{\beta,\mathrm{tag}}(t) \right].
$$
\end{cor}
Note that the term $b_{\beta}$ is negative when $\beta<-3/2$ and null for $\beta=-3/2$. When $\beta>-3/2$ we cannot use directly Lemma \ref{prop:aenuinfinite}. In these cases the $\gamma_1$ in the asymptotic expansion of $\phi$ (with the notation of Lemma \ref{prop:aenuinfinite}) is equal to $-1-\beta$ and therefore is smaller than $1/2$: hence there will be additional terms in or in front of the exponential (more and more terms as $\beta \uparrow -1)$.

\bigskip

\bibliographystyle{acm}
\bibliography{frag}

\end{document}